\documentclass[12pt]{amsart}
\input{diagrams.tex}
\diagramstyle[scriptlabels,height=8mm,width=8mm]
\newcommand{\query}[1]%
{\mbox{}\marginpar{\raggedright\hspace{0pt}{\small\em #1}}}%
\theoremstyle{plain}

\newtheorem{thm}{Theorem}[section]
\newtheorem{cor}[thm]{Corollary}
\newtheorem{lem}[thm]{Lemma}
\newtheorem{prop}[thm]{Proposition}

\theoremstyle{definition}
\newtheorem{defi}[thm]{Definition}
\newtheorem{conj}[thm]{Conjecture}
\newtheorem{conv}[thm]{Convention}
\newtheorem{nota}[thm]{Notation}
\newtheorem{rem}[thm]{Remark}
\newtheorem{exa}[thm]{Example}
\newtheorem{sit}[thm]{}
\newtheorem{que}[thm]{Question}

\newcommand{\brem}{\begin{rem}}
\newcommand{\erem}{\end{rem}}
\newcommand{\bexa}{\begin{exa}}
\newcommand{\eexa}{\end{exa}}
\newcommand{\bdefi}{\begin{defi}}
\newcommand{\edefi}{\end{defi}}
\newcommand{\bcor}{\begin{cor}}
\newcommand{\ecor}{\end{cor}}
\newcommand{\blem}{\begin{lem}}
\newcommand{\elem}{\end{lem}}
\newcommand{\bconv}{\begin{conv}}
\newcommand{\econv}{\end{conv}}
\newcommand{\bconj}{\begin{conj}}
\newcommand{\econj}{\end{conj}}
\newcommand{\bprop}{\begin{prop}}
\newcommand{\eprop}{\end{prop}}
\newcommand{\bthm}{\begin{thm}}
\newcommand{\ethm}{\end{thm}}
\newcommand{\bnota}{\begin{nota}}
\newcommand{\enota}{\end{nota}}
\newcommand{\bsit}{\begin{sit}}
\newcommand{\esit}{\end{sit}}

\newcommand{\orbit}{{\operatorname{O}}}

\newcommand{\Sing}{\operatorname{Sing}}
\newcommand{\Spec}{\operatorname{Spec}}
\newcommand{\Frac}{\operatorname{Frac}}
\newcommand{\red}{{\operatorname{red}}}
\newcommand{\Proj}{\operatorname{Proj}}
\newcommand{\ord}{\operatorname{ord}}
\newcommand{\id}{\operatorname{id}}
\newcommand{\Der}{\operatorname{Der}}

\def\fC{{\mathfrak C}}
\def\fH{{\mathfrak H}}

\def\cK{{\mathcal K}}
\def\cL{{\mathcal L}}
\def\cO{{\mathcal O}}
\def\cX{{\mathcal X}}

\newcommand{\pP}{{\mathbb P}}

\newcommand{\C}{{\mathbb C}}
\newcommand{\Q}{{\mathbb Q}}
\newcommand{\Z}{{\mathbb Z}}
\newcommand{\N}{{\mathbb N}}

\newcommand{\rmap}{\mbox{\hskip4pt-\hskip1pt-\hskip1pt-%
\hskip1pt-}\rhla\hskip2pt}

\def\lto{\longrightarrow}
\def\hto{\hookrightarrow}

\def\CM{Cohen-Macaulay}
\newcommand{\la}{\label}

\newcommand{\G}{{\Gamma}}

\newcommand{\p}{{\partial}}

\newcommand{\be}{\begin{eqnarray}}
\newcommand{\ee}{\end{eqnarray}}
\newcommand{\no}{\noindent}

\edef\qqed{\qed}
\def\qed{\qqed \medskip}

\title{Rational curves and rational
singularities}

\author{Hubert Flenner}
\address{Fakult\"at f\"ur Mathematik,
Ruhr Universit\"at Bochum,
Geb.\ NA 2/72,
Universit\"ats\-str.\ 150,
44780 Bochum, Germany}
\email{Hubert.Flenner@ruhr-uni-bochum.de}

\author{Mikhail Zaidenberg}
\address{Universit\'e
Grenoble I, Institut Fourier, UMR 5582 CNRS-UJF, BP 74,
38402 St.\ Martin
d'H\`eres c\'edex, France}
\email{zaidenbe@ujf-grenoble.fr}

\thanks{
\mbox{\hspace{11pt}}{\it 1991 Mathematics Subject Classification}:
14J17, 14L30,  13H10.\\
\mbox{\hspace{11pt}}{\it Key words}: $\C^*$-action, $\C_+$-action,
graded algebra, rational curve, polynomial curve, affine surface, rational
singularity, quotient singularity}

\date{\today}

\begin{document}

\begin{abstract}
We study rational curves on  algebraic varieties, especially
on normal affine varieties  endowed with a
$\C^*$-action. For varieties with an isolated singularity, we show
that the presence of sufficiently many rational curves
outside the singular point strongly affects
the character of the singularity.
This provides an explanation of classical results due to
H.  A. Schwartz and G. H. Halphen
on polynomial solutions of the generalized Fermat equation.
\end{abstract}

\maketitle
{\footnotesize \tableofcontents}

\section*{Introduction}

In 1873 H. A. Schwartz \cite{Schw}
(in his studies on the Gauss hypergeometric equation)
found polynomial solutions in coprime polynomials $x(s),\,y(s),\,z(s)$
of the generalized Fermat equation
\be\label{plat} x^p+y^q+z^r=0\,\ee
for every Platonic triple $(p,q,r)$ with $p,q,r\ge 2,\,\,\,1/p+1/q+1/r>1$.
These solutions are given
(up to constant factors, in the notation as in \cite{Schw}) by the
following identities.

\no For a dihedral triple $(p,q,r)=(2,2,d)$:
$$(s^d+1)^2-(s^d-1)^2=4s^d\,.$$
For the tetrahedral triple $(p,q,r)=(2,3,3)$:
$$12\sqrt{3}(s(1+s^4))^2=(1+2\sqrt{3}s^2-s^4)^3-(1-2\sqrt{3}s^2-s^4)^3\,.$$
For the octahedral triple $(p,q,r)=(2,3,4)$:
$$(1-33s^4-33s^8+s^{12})^2=(1+14s^4+s^8)^3-
4\cdot s^3(s(1-s^4))^4\,.$$
For the icosahedral triple $(p,q,r)=(2,3,5)$:
$$[\varphi_{30}(s)]^2=[\varphi_{20}(s)]^3-4^3\cdot 3^3\cdot
[\varphi_{12}(s)]^5\,,$$
where
$$
\begin{array}{ccc}
\varphi_{12}(s) & = & s(1-11s^5-s^{10})\\
\varphi_{20}(s) & = & 1+228s^5+494s^{10}-228s^{15}+s^{20}\\
\varphi_{30}(s) & = & 1-522s^5-10005s^{10}-10005s^{20}+522s^{25}+s^{30}\,.
\end{array}
$$
Schwartz' solutions are by no means unique, some others having been known
to Euler (1756), Hoppe (1859), Liouville (1879), etc. \cite{DarGr}.
On the other hand, in 1880 Halphen \cite{Ha} showed that
the generalized Fermat equation (\ref{plat}) has no solution
in non-constant coprime polynomials when $1/p+1/q+1/r\le 1$.
We highly recommend to consult \cite{Kl} for a historical account
on the subject, \cite{BalDw, Beu} for a modern one,
and \cite{DarGr} for both. (Notice that the above
solutions are given there via homogeneous forms in two variables.)

If $(a,b,c)\in\C^3$ is a nonzero constant solution of (\ref{plat})
and $f\in \C[s]$
then
clearly $(af^{M/p},\,bf^{M/q},\,cf^{M/r})$ with $M:={\rm lcm}(p,q,r)$ is
a solution of (\ref{plat}) which we call {\em trivial}.
The existence of non-trivial solutions is equivalent to
the quasirationality of the singularity of the Pham-Brieskorn
surface $V_{p,q,r}$ given in
$\C^3$ by the equation (\ref{plat})
\cite[L. 7]{KalZa}. The latter can be expressed explicitly in terms of the
exponents $(p,q,r)$
as in (a) below, whereas the classical results described above
admit a geometric interpretation as in (b).

\bthm\label{PHBR} $\,$
\begin{enumerate}
\item[(a)] \cite[Cor.\ of L.~8]{BarKau, Ev, KalZa}
There exists a polynomial curve $C\subseteq V_{p,q,r}$
(i.e.\ $C=f(\C)$, where $f:\C\to
V_{p,q,r}$ is a non-constant morphism) that is not
contained in an orbit closure of the natural
$\C^*$-action on\footnote{This corresponds to a
non-trivial polynomial solution of (\ref{plat}).}
$V_{p,q,r}$ if and only if either
\begin{enumerate}
\item[(i)] one of the exponents $p,\,q,\,r$
is coprime with the other two, or
\item[(ii)] gcd$(p,q,r)=2$ and $p/2,\,q/2,\,r/2$
are coprime in pairs.
\end{enumerate}
\smallskip

\item[(b)] {\em (Schwartz-Halphen)} There exists
a polynomial curve $C$ in $V_{p,q,r}$
not passing through the origin if
and only if $(p,q,r)$
is a Platonic triple (or, what is equivalent,
iff the singularity of the surface $V_{p,q,r}$ at the
origin is a rational double point).
\end{enumerate}
\ethm

In sections 1 and 2 we obtain similar results
for more general
quasihomogeneous
affine algebraic varieties, i.e., varieties endowed
with an effective
$\C^*$-action\footnote{To avoid possible confusion,
notice that algebraic group theorists use this same term in
a different sense. Namely, they call a variety
{\it quasihomogeneous} if it
possesses an algebraic group action with a
Zariski open dense orbit.}.
The most complete results concerns surfaces with many
rational curves. By a {\em rational curve}
we always mean a curve with function field isomorphic to
$\C(t)$ or, in other words, the normalization $C'$ of $C$
admits an open embedding into $\pP^1$. Our results can be
formulated as follows (see Theorems
\ref{newthm0},
\ref{mainthm} and
\ref{cycl}).

\bthm\label{mt}
Let $V$ be a normal affine algebraic
surface
over $\C$ with a good $\C^*$-action\footnote{A regular
$\C^*$-action on an algebraic variety is called {\em good}
if it possesses a unique attractive fixed point.}
and with an isolated
singularity $p\in V$.
Then the following  hold.

\begin{enumerate}
\item[(a)] If there exists a closed rational
curve $C$ on $V$
which does not pass through the singular point $p$
then $(V,p)$ is a rational singularity.

\item[(b)] There exists a polynomial curve
$f:\C\to V$ on $V$
which does not pass through the point $p$
if and only if $(V,p)$ is a quotient singularity.

\item[(c)] The surface $V$ admits a nontrivial
regular $\C_+$-action
if and only if $(V,p)$ is a cyclic quotient singularity.
\end{enumerate} \ethm

The implication "$\Rightarrow$" in (c) is a special
case of a more general result of Miyanishi \cite{Miy}.
We do not know whether the converse to (a) is also true.
Notice that if a surface singularity $(V,p)$ as in the theorem
is Gorenstein (for instance,
if this is an isolated complete intersection surface singularity)
then conditions (a) and (b) above are equivalent (Theorem
\ref{lines}).

More generally, (a) is true in any dimension and
for not necessarily quasihomogeneous varieties.
Namely, we have the following result (see Theorem \ref{newthm0}).

\bthm\label{mth} If $V$ is a
scheme (or a Moishezon variety) with an isolated Cohen-Macaulay singularity
$p\in V$ such that

\begin{enumerate}
\item[(*)] a Zariski open subset of $V$
can be covered by closed rational curves $C$
which do not pass through $p$
\end{enumerate}
then
$(V,p)$ is a rational singularity. \ethm

Under the assumption (*) we say that $V$ is {\em uniruled off $p$}.
Notice that the assumption $p\not\in C$
is essential as is easily seen by the example of non-rational
quasihomogeneous isolated hypersurface singularities.

As an application, consider the
hypersurface $V$ in $\C^n$  given by an equation
$$P:=uv-p(x_1,\ldots,x_{n-2})=0.$$ It admits  an
effective
$\C_+$-action induced by the locally nilpotent derivation
$\partial$ on $\C[x_1,\ldots,x_n]$ with
$$
\partial(u):=0\,,\quad
\partial (v):=\frac{\partial p}{\partial x_1}\,,\quad
\partial(x_1):=u\quad\mbox{and}\quad
\partial(x_i):=0\mbox{ for }i\ge 2.
$$
Assume that $V$ has an isolated singularity at $0$.
The general orbits of the $\C_+$-action constitute a family of polynomial
curves on $V$
not passing through $0$. Thus we recover the result of Viehweg
\cite{Vie} saying that $(V,0)$ is rational.

Throughout the paper we use log-canonical $L^2$-forms
on singular complex spaces as considered in \cite{FlZa1};
we recall their definition and some useful facts
in subsection 1.1. As for the classification results
\ref{glob}, \ref{cycl} and \ref{cyquo}
concerning "good" quasihomogeneous
surfaces  with a $\C_+$-action
(as well as for more general ones)
see also the forthcoming paper \cite{FlZa2}.

\section{Rational curves on quasihomogeneous varieties}

\subsection{Preliminaries: log-canonical sheaves}
Here we recollect some notions and facts from \cite{FlZa1}
that we need in the sequel.

\bsit \label{cf}
Let $X$ be a normal complex space, and let
$\sigma : Y \to X$ be a resolution of singularities
such that $E:=\sigma^{-1}(\Sing X)_\red$ is a divisor
with simple normal crossings (SNC divisor in brief). Define
$${\cL}_{Y,E}^{2,m}:={\cO}_Y(mK_Y+(m-1)E)\qquad {\rm and}\qquad
  {\cL}_{X}^{2,m}:=\sigma_* {\cL}_{Y,E}^{2,m}\,;$$
similarly,
$${\cL}_{Y,E}^{m}:={\cO}_Y(mK_Y+mE)\qquad {\rm and}\qquad
  {\cL}_{X}^{m}:=\sigma_* {\cL}_{Y,E}^{m}\,.$$
As explained in \cite{Sak},
the forms in $H^0(Y,\cL_{Y,E}^{2,m})$ are just the meromorphic
$m$-canonical forms on
$Y\backslash E$ which are locally $L^2$
at the points of $E$.
The sheaf ${\cL}_{X}^{2,m}$ does not
depend on the choice of resolution
\cite[Sect.\ 1.1]{FlZa1}, and
its sections
over an open subset $U\subseteq X$ are the sections in
$H^0(U_{\rm reg},\,\omega_X^{\otimes m})$ which are locally
$L^2$ on $U$ \cite{La2, Bur, KiWa}.
In general, ${\cL}_{X}^{2,m}\subseteq \cO_X(mK_X)$.
Assuming that $K_X$ is a $\Q$-Cartier divisor,
the equality holds
if and only if $X$ has at most
log-terminal singularities \cite[Prop.\ 1.17(c)]{FlZa1}.
Recall that in dimension 2
log-terminal singularities
are just quotient singularities
(Kawamata \cite{Kaw}; see e.g., \cite[3.0.1]{FA}).

If $D\subseteq X$ is a closed reduced subspace then we may also
consider the log-canonical sheaves ${\cL}_{X, D}^{m}$
and
${\cL}_{X,D}^{2,m}$ taking
$E:=\sigma^{-1}(\Sing X\cup D)_\red$ in the above definitions.
\esit

Of particular importance is the sheaf $\cK_X:=\cL_X^{2,1}=\sigma_*
\omega_{X'}\subseteq \omega_X$
of $L^2$-canonical forms
which was called in \cite{GrRi} the {\em canonical
sheaf}\/ of $X$. A basic result is the
Grauert-Riemenschneider vanishing theorem along with its
generalizations due to Koll\'ar and Moriwaki \cite{Moriw}
(see also M. Saito
\cite{Sai} and Arapura \cite{Ara}),
which we formulate for later purposes
as follows.

\bprop\label{Kollar}
{\em (a)  (\cite{GrRi}, \cite[Thm. 3.2]{Moriw})}
Let $\pi:X'\to X$ be a birational proper morphism of normal
complex spaces. Then $\pi_*(\cK_{X'})\cong \cK_X$ and
$R^i\pi_*(\cK_{X'})=0$ for $i\ge 1$.

{\rm (b)  (\cite[Thm. 2.1(i)]{Kol1},
\cite[Thm. 3.2]{Moriw})}
Let $\pi : X'\to X$
be a morphism
of irreducible normal complex spaces which admits a factorization
$$
\pi: X'\stackrel{p}{\lto} Y\stackrel{q}{\lto} X
$$
such that

(i) $p$ is  projective, i.e.\ there is a $p$-ample
line bundle $\cO_{X'}(1)$;

(ii) $q$ is a Moishezon morphism.\footnote{That is, for any point
$x\in X$ there exist a neighborhood $U$ of $x$ in $X$,
a compact Moishezon variety $V$ and a closed embedding
$\varphi :q^{-1}(U)\hookrightarrow U\times V$ such that
$q\big | q^{-1}(U) = {\rm pr}_1\circ\varphi$, see
\cite[Sect.\ 3]{Moriw}.}

\no Then all direct image sheaves
$R^i\pi_*(\cK_{X'})\,\,\,(i\ge 0)$
on $X$ are torsion free.
\eprop

\proof
(a) In the case that $\pi$ is projective and $X'$ is smooth,
Moriwaki's result \cite[Thm.\ 3.2]{Moriw} gives that the sheaf
$R^i\pi_*(\omega_{X'})$, $i\ge 1$, is torsion free on $X$ and
so, being concentrated on $\Sing X$, it must vanish. In the
general case we can find a desingularization $\sigma:X''\to
X'$ such that $\pi\circ\sigma$ is a projective morphism.
As $\sigma$ is then projective as well, by the first part of
the proof the sheaves $R^i\sigma_*(\omega_{X''})$ and
$R^i(\pi\circ\sigma)_*(\omega_{X''})$ vanish for
$i\ge 1$. By definition, $\sigma_*(\omega_{X''})=\cK_{X'}$
and so $R^i\pi_*(\cK_{X'})\cong
R^i(\pi\circ\sigma)_*(\omega_{X''})=0$ for $i\ge 1$, proving (a).

(b) If $\pi$ is a Moishezon morphism and $X'$ is smooth then (b) is
shown in \cite[3.2]{Moriw}. The general case is reduced to this as
follows: replacing $Y$ by the image of
$p$, we may suppose that $p$ is surjective and $Y$ is irreducible.
According to
$[loc.cit$, 3.3], locally with respect to $X$ we can find an
irreducible normal complex space $\tilde Y$ and a proper
bimeromorphic morphism $\mu:\tilde Y\to Y$ such that the composition
$\tilde q:=q\circ\mu$ is projective. Let $\tilde X'\to \tilde
Y\times_Y X'$ be a projective morphism that desingularizes the
fiber product. Then all maps in the diagram
\begin{diagram}
\tilde X' & \rTo^{\tilde p} & \tilde Y\\
\dTo<{\mu'} && \dTo<{\mu} &\rdTo^{\tilde q}\\
X'& \rTo^p & Y & \rTo^q & X
\end{diagram}
except possibly $q$ are projective. By $loc.cit.$\
the  direct image sheaves $R^i(\tilde q\circ\tilde
p)_*(\omega_{\tilde X'})$ are torsion free for all $i$. By
construction, the map $\mu'$ is proper and bimeromorphic and so,
using (a), $R\mu'_*(\omega_{\tilde X'})\cong
\cK_{X'}$ in the derived category. Thus
$$
R^i(\tilde q\circ\tilde p)_*(\omega_{\tilde X'})\cong
R^i\pi_*(\cK_{X'})
$$
for all $i$, which gives (b).
\qed

\bsit\label{pgr}
Let $(X,\,p)$ be an isolated singularity of a normal complex space.
Following \cite{KiWa} we define the
{\it $m$-th $L^2$-plurigenus} of $(X,\,p)$ to be
$$\delta_m\,(X,\,p) = \dim_{\C}\,
\left[\omega_{X,\,p}^{[m]}/({\mathcal
L}_{X}^{2,m})_p\right]\,,$$
where $\omega_X^{[m]}$ denotes the reflexive hull of
$\omega_X^{\otimes m}$.
\esit

We make frequent use of the following results.

\bprop\label{Kempf}
{\rm (Kempf \cite[Prop.\ on p.\ 50]{KKMSD})}
Let $(X,p)$ be a singularity of a normal complex space.
Then $(X,p)$ is
a rational singularity if and only if it is Cohen-Macaulay
and $\cK_{X,p}\cong \omega_{X,p}$.
In the case of an isolated normal singularity $(X,p)$
the latter condition is equivalent to $\delta_1\,(X,\,p) =0$.
\eprop

Note that in $loc.cit.$\ this proposition is shown only for
algebraic singularities. However, as remarked in
\cite[1.17(a)]{FlZa1} the result holds as well in the complex
analytic category.

\bprop\label{Watanabe}
{\rm (Ki.\ Watanabe \cite[Thm.\ 3.9]{KiWa})}
If $(X,p)$ is an isolated normal surface singularity
then $\,\delta_m\,(X,\,p) =0\,\,\forall m\ge 1$ if and only if
$(X, p)$ is a quotient  singularity.
\eprop

\bsit\label{plge} Let   $V$ be a
smooth variety  with a
smooth compactification  $\bar V$  by an SNC-divisor
$D=\bar V\backslash V$. The {\it logarithmic
plurigenera}
$\bar p_m(V)$ are defined by
$$\bar p_m(V)=\dim\,H^0(\bar V,\,{\mathcal
O}_{\bar V}(mK_{\bar V}+mD))\,,\,\,\,m\ge 1,$$
whereas the {\it logarithmic Kodaira dimension} $\bar k\,(V)$ of $V$ is
$$
\bar k\,(V)=
\left\{ \begin{array}{cc}
-\infty & {\rm if}\quad \bar p_m(V)=0\quad\forall
m\in\N\\
\min\,\{k\in \N : \limsup\limits_{m\to\infty}\bar
p_m(V)/m^k<\infty\} &
\quad{\rm otherwise.}
\end{array} \right .
$$
For a possibly singular variety $V$ one defines $\bar
p_m(V):=\bar p_m(V')$ and
$\bar k\,(V):=\bar k\,( V')$,
where $ V' \to V$ is a  resolution of singularities.
Let $\bar V$ be a compactification of $V$ and $D:=\bar
V\backslash V$. If $\pi:\bar V'\to \bar V$ is a resolution of
singularities such that $D':=\pi^{-1}(D)_\red$ is an SNC
divisor then by our definitions there are inclusions
$$
H^0(\bar V',\,\cO_{\bar V'}(m(K_{\bar V'}+D)))
\subseteq
H^0(\bar V,\,\cL^m_{\bar V,D})
\subseteq
H^0(V,\,\cL^m_V)
$$
(cf.\ \ref{cf}). Thus, if the group on the right or in the
middle vanishes for $m\ge 1$ then $\bar p_m(V)=0$ and, in
particular, $\bar k(V)=-\infty$.
\esit

In the case of quasihomogeneous singularities
we have the following facts.

\bthm\label{dell}  Let
$A=\bigoplus_{\nu\ge 0} A_{\nu}$ be a normal $\C$-algebra of
finite type and assume that the corresponding
affine variety $V=\Spec A$
has at most isolated singularities. Then the following hold.

{\rm (a)  (\cite[2.22 (b)]{FlZa1})}
$\delta_m (V,p)=\dim_{\C}\,(\omega_A^{[m]})_{\le
0}$, where $\omega_A=H^0(V,\cO(K_V))$
is the dualizing module of $A$.
In particular, if $V$ is Cohen-Macaulay
then it has at most rational singularities
if and only if
$(\omega_A)_{\le 0}=0$.

{\rm (b)  (\cite[2.22 (a)]{FlZa1})} If
$A_0\ne \C$ (in particular, if $V$
has at least two singular points)
then  $\delta_m (V,p)=0$
for all $m\ge 1$ and $p\in\Sing V$.

{\rm (c) (\cite[2.26 (a), (b)]{FlZa1})}
If $A_0=\C$ then
$\bar p_m(V\backslash \{p\})=
\dim\,(\omega_A^{[m]})_{0}\,$, where $\{p\}:=V(A_+)$
with $A_+:=\bigoplus_{\nu> 0} A_{\nu}$. Consequently,
$\bar k(V\backslash\{p\})=-\infty$ if and only if
$(\omega_A^{[m]})_0=0$ for all $m\ge 1$ or, equivalently,
if and only if $\delta_m(V,\,p)=0$ for all $m\ge 1$.
\ethm

\bsit\label{ci}
Let $V=\Spec A$ be a complete intersection
of  dimension
$n$ given in $\C^{n+s}$ by polynomials
$p_1,\dots, p_{s}\in \C^{n+s}$
that are  quasihomogeneous of degrees
$d_1,\dots, d_s$ with respect to
weights $w_j> 0\,\,(j=1,\dots,n+s)$ i.e.,
$$
p_i(\lambda^{w_{1}}x_1,\dots,\lambda^{w_{n+s}}x_{n+s})=
\lambda^{d_i}p_i(x_1,\dots,x_{n+s}),\qquad i=1,\dots,s\, .
$$
It is well known that $\omega_A=A[N_A]$ (i.e.,
$(\omega_A)_{\nu}=A_{\nu+N_A}$) and
$(\omega_A^{[m]})_0 =A_{mN_A}$, where
$$
N_A (=N_V):=\sum_{i=1}^{s} d_i -
\sum_{j=1}^{n+s}w_{j}\,,
$$
see e.g.\  \cite[2.2.8 and 2.2.10]{Go} or \cite[p.\ 42]{Fl}
for a simple argument.
Applying the foregoing results to the case of such complete
intersections we get the following corollary.
\esit

\bcor\label{cor ci}
  {\rm ( \cite{Fl}, \cite{Mora},
\cite[2.23, 2.26(c), 2.28]{FlZa1})}
Suppose that $V$ has an isolated
singularity at the origin
$p:={\bar 0}\in \C^{n+s}$. Then
$$ \delta_m (V,p)=\dim_{\C}\,\sum_{\nu\le 0} A_{\nu+mN_A},\qquad
{\bar p}_m(V\backslash \{p\})=\dim_{\C}\,A_{mN_A}\,,$$
and so
$$
\bar k(V\backslash \{p\})=
\left\{ \begin{array}{ccc}
-\infty & {\rm iff}\quad N_A<0\\
0 & {\rm iff}\quad N_A=0\\
\dim V -1 & {\rm iff}\quad N_A>0\,.
\end{array} \right .
$$ Consequently, $(V,p)$ is a rational singularity if and only if
$N_A<0$.
\ecor

\subsection{Canonical forms and rational curves} It is well
known that complete smooth uniruled varieties have no
$m$-canonical forms (see e.g.\ \cite[IV, Cor. 1.11]{Kol}). In
the following proposition (b) and (c) we show a similar
result for varieties that admit ``many'' affine rational
curves. Notice that the assumption in (a) below is just that
the variety is uniruled (cf.\
\cite[IV, Def. 1.1]{Kol}).

\bprop\label{forms}
Let $\bar X$ be a complete normal algebraic variety
over $\C$, let $D\subseteq X$ be a closed reduced subspace and
assume that $X:=\bar X\backslash D$ is smooth.
  \begin{enumerate}
\item[(a)] If for a general point $x$ of $X$ there is a non-constant
rational  curve $f: \pP^1\to X$ passing through $x$ then $H^0(X,
\omega^{\otimes m}_{X})=0$ for all $m\ge 1$. In particular,
$\bar k(X)=-\infty$.
\item[(b)]  If for a general point $x$ of $X$ there is a non-constant
curve $f:\C\to X$ passing through $x$ then
$H^0(\bar X, \cL^{m}_{\bar
X, D})=0$ for all $m\ge 1$. In particular, $\bar k(X)=-\infty$.
\item[(c)]  If for a general point $x$ of $X$ there is a non-constant
curve $f:\C^*\to X$ passing through $x$ then
$H^0(\bar X, \cL^{2,m}_{\bar
X, D})=0$ for all $m\ge 1$.
\end{enumerate}
\eprop

\proof
The second assertions in (a) and (b) follow from \ref{plge}. For the
remaining statements, consider the Hilbert scheme $\fH$ of $\bar
X$.  As
$\fH$ has a countable number of
irreducible components there is a closed subvariety $S'$ of $\fH$ with
the following properties.

\begin{enumerate}
\item[(i)]  The general point of $S'$ is represented by a complete
rational curve $\bar C$ in $\bar X$ such that the normalization of $C:=
\bar C\backslash D$ is isomorphic either to $\pP^1$ (or to $\C$
or $\C^*$ in cases (b) and (c), respectively).
\item[(ii)]  The curves parametrized by $S'$ cover $\bar X$.
\end{enumerate}
After passing to a suitable subvariety  of $S'$ we may also
assume that $\dim S'=\dim X- 1$.
Consider the incidence subspace
$$
\bar\cX':=\{(x,[\bar C])\in \bar X\times S'\, |\, x\in \bar C\}
$$
of $\bar X\times S'$ which is flat over $S'$. Passing to the
normalizations $\bar\cX$ of $\bar\cX'$ and $S$ of $S'$ we get morphisms
\begin{diagram}[h=7mm]
\bar\cX & \rTo^{f} & \bar X\\
\dTo<{\pi} \\
S
\end{diagram}
By construction the fibers of $\pi$ are all curves and its
general fibers  are isomorphic to $\pP^1$. Moreover,
$f$ is dominant and $\dim \bar\cX=\dim \bar X$.

The subspace $\Delta:=f^{-1}(D)_{\rm
red}$ of $\bar\cX$ satisfies $\deg
\Delta / S \le 0$ (resp., $\le 1$ and $ \le 2$ in case (b) and
(c)). Using the functoriality of $m$- and $(2,m)$-canonical
forms (see \cite[Sect.\ 1]{FlZa1}), with
$\cX:=\bar\cX\backslash \Delta$ we have inclusions
$$
f^*(\omega^{m}_X)\subseteq \cL^{m}_{\cX}\,\,,
\quad
f^*(\cL^{m}_{\bar X, D})\subseteq \cL^{m}_{\bar
\cX, \Delta}\quad\mbox{and}\quad
f^*(\cL^{2,m}_{\bar X, D})\subseteq \cL^{2,m}_{\bar
\cX, \Delta}\,.
$$
Thus it is sufficient to show that for $m\ge 1$
$$
H^0(\cX, \cL^{m}_{\cX})=0\,\,,\quad
H^0(\bar \cX, \cL^{m}_{\bar
\cX, \Delta})=0\quad
\mbox{and}\quad
H^0(\bar \cX, \cL^{2,m}_{\bar
\cX, \Delta})=0
$$
in case (a), (b) and (c), respectively.
In case (a) this follows easily from the fact that a general fiber
$F\simeq \pP^1$ of $\pi$ does not intersect $\Delta$, whence
$$
\deg \cL^{m}_\cX|F=\deg \omega^{\otimes m}_\cX|F=-2m<0.
$$
Similarly, in case (b) we have
$$
\deg \cL^{m}_{\bar\cX, \Delta}|F=m(-2+\deg \Delta / S)<0,
$$
whereas in the situation of (c)
$$
\deg \cL^{2,m}_{\bar\cX,\Delta}|F=-2m+(m-1)\deg \Delta/S<0.
$$
\qed

\brem\label{alldir}
Let $\bar X$ and $S$ be algebraic varieties over $\C$ and let $\bar
f_s:\bar C_s\to\bar X$ be a family of complete curves in $\bar X$
parameterized by $S$ that cover $\bar X$ in the sense that there is a
proper smooth family $\pi: \bar\cX\to S$ of curves and a
dominant morphism
$\bar f:\bar\cX\to\bar X$ such that for $s\in S$ we have $\bar
C_s=\pi^{-1}(s)$ and $f|\bar C_s=\bar f_s$.

Let $D\subseteq \bar X$ be a closed reduced subspace and let $X:=\bar
X\backslash D$. With $\Delta:=\bar f^{-1}(D)_\red$
consider the family of curves
$$
f_s:C_s:=\bar C_s\backslash \Delta\to X
$$
with $s\in S$, where $f_s:=\bar f_s|C_s$. We call
$\kappa:=\deg \Delta/S$ the {\em number of punctures } of the family
$f_s$.

If $\bar X$ is smooth and $D$ is a divisor then for general $s\in S$
we have estimates
\be
(K_{\bar X}+D)\bar f_{s*}(\bar C_s)\le 2g-2+\kappa\,,\\
K_{\bar X}\bar f_{s*}(\bar C_s)\le 2g-2\,,
\ee
where $g$ is the genus of $\bar C_s$.
\erem

\proof
We may assume that $S$ is smooth of dimension $\dim
X-1$. As above
$\bar f^*(\omega_{\bar X}(D))\subseteq
\omega_{\bar\cX}(\Delta)$. For general $s\in S$ this restricts
to a non-zero map $\bar f_s^*(\omega_{\bar X}(D))\to
\omega_{\bar C_s}(\Delta\cap \bar C_s)$; note that
$\omega_{\bar\cX}(\Delta)|\bar C_s\cong \omega_{\bar C_s}(\Delta\cap
\bar C_s)$ as the normal bundle of $\bar C_s$ in $\bar\cX$ is trivial.
Hence
$$
(K_{\bar X}+D)\bar f_{s*}(\bar C_s)=\deg \bar f_s^*(\omega_{\bar
X}(D))\le \deg \omega_{\bar C_s}(\Delta\cap\bar C_s)=2g-2+\kappa,
$$
proving the first inequality. Applying this to the case $D=\emptyset,\
\Delta=\emptyset$ gives the second inequality.
\qed

\subsection{Rational curves and rational singularities}
In this subsection we
prove theorem \ref{mth} from the introduction.

\bthm\label{newthm0} Let $p\in X$ be an
isolated Cohen-Macaulay
singularity of an algebraic scheme
(resp., of a compact Moishezon variety) $X$.
If for a general point $x\in X$ there exists
a closed rational curve $C\subseteq X$
passing through $x$ but not through $p$ then
$(X,p)$ is a rational
singularity.\ethm

\proof
  In the case of an algebraic scheme
we may assume that $X$ and $C$ are projective
(first passing to an affine neighbourhood of
$p$ in $X$ and then to its
normalized projective closure).
Thus it is enough to consider the case where $X$ is
a compact Moishezon variety.
Let $\fH$ be the Hilbert scheme of $X$. As
in the proof of \ref{forms} there is an irreducible compact subvariety
$S'\subseteq \fH$ of dimension $\dim S'=\dim X-1$ such that the following
hold.

(a) A general point of $S'$ represents a rational curve not passing
through $p$;

(b) the curves in $S'$ cover $X$.

\noindent As above let
$\cX':=\{(x, [D])\in X\times S'\,\big |\,x\in D\}$
be the incidence subvariety of $X\times S'$ which is
flat over $S'$. Passing to the normalizations
$\cX$ of $\cX'$ and $S$ of $S'$ we get morphisms
\begin{diagram}[h=7mm]
\cX & \rTo^{f} & X\\
\dTo<{\pi} \\
S
\end{diagram}
By construction the fibers of $\pi$ are all curves,
its general fibers are isomorphic to $\pP^1$ and
$f$ is finite on every fiber of $\pi$.
As general fibers of $\pi$ do not pass through
the singular point $p\in X$ and the preimage
$A:=f^{-1}(p)$ does not contain a whole fiber of $\pi$,
this analytic subset $A$
has codimension $\ge 2$ in $\cX$. Moreover  $\pi : A\to S$
is a finite morphism and $\cX\to X\times S$
is finite.

Consider now the sheaves
$$
\cK_{\cX}:=\cL_{\cX}^{2,1},\quad \cK_{X}:=\cL_{X}^{2,1}
\quad {\rm and}\quad \cL:=f^*(\omega_X)\cap
\omega_{\cX}\subseteq
\omega_{\cX}\,.
$$
Let us show that the quotient $\omega_X/\cK_X$
vanishes in a neighborhood of $p\in X$.
Indeed,  $\cL_x=f^*(\omega_X)_x$ if $f(x)$ is a regular
point of $X$. As codim$\,A\ge 2$ the sheaf
$\cL$ coincides with $f^*\omega_X$ in a neighborhood of $A$.
Moreover, by the functoriality
of $L^{2,m}$-canonical forms \cite[Sect.\ 1]{FlZa1}
$f^*\cK_X\subseteq \cK_{\cX}$.
Thus the sheaf $\fC$ in the exact sequence
\be\label{exse}
0\to \cK_{\cX}\to\cK_{\cX}+\cL\to \fC:=
(\cK_{\cX}+\cL)/\cK_{\cX}\cong
\cL/(\cK_{\cX}\cap\cL)\to 0\ee
is a subquotient of $f^*\omega_X/f^*\cK_X$.
Since $p\in X$ is an isolated singularity,  and so
$\cK_X=\omega_X$ in $X\backslash \{p\}$ near $p$,
it follows that
the support of $\fC$
is contained in
$A\cup A'$ where $A'\subseteq \cX$ is a
closed analytic subset disjoint from $A$.
In particular, $\fC\cong\fC|A\oplus\fC|{A'}$.
The sheaf
$\pi_*(\cK_{\cX}+\cL)$ is contained in
$\pi_*(\omega_{\cX})$, whence it vanishes
as the general
fibers of $\pi$ are isomorphic to $\pP^1$.
Taking $\pi_*$ of the sequence (\ref{exse})
we get that $\pi_*(\fC)$ is a subsheaf of
$R^1\pi_*(\cK_{\cX})$.

We claim that $\pi$ can be written as the composition of a
projective morphism and a  Moishezon morphism.  Indeed, as $X$ is
Moishezon and $\cX'$ is contained in $X\times S'$, the natural
projection $\pi' :\cX'\to S'$ and then by base change the projection
$\cX'\times_{S'}S\to S$ are Moishezon as well.
As $\cX\to \cX'\times_{S'}S$ is a finite map it is in particular
projective. Thus
$$
\pi: \cX\lto \cX'\times_{S'}S\lto S
$$
is the desired decomposition.

Applying Proposition \ref{Kollar} (b) the  sheaf
$R^1\pi_*(\cK_{\cX})$ on $S$ is torsion free and so its
subsheaf
$\pi_*(\fC)$ is torsion free as well.
Hence $\fC|A=0$ since otherwise
$\pi_*(\fC|A)$ would be a torsion subsheaf of
$\pi_*(\fC)\cong\pi_*(\fC|A)\oplus\pi_*(\fC|{A'})$.
This shows that
$f^*(\omega_X)= \cL\subseteq\cK_{\cX}$
in a neighborhood of $A$. Applying \cite[Prop.\
1.9 (b))]{FlZa1} we obtain that $\omega_X=\cK_{X}$
in a neighborhood of $p\in X$. Hence by
\ref{Kempf} the singularity $(X,p)$
is rational,
as stated. \qed

\bcor\label{agp}
Let $p\in X$ be an isolated \CM\ singularity of an affine
algebraic variety. If $X$ admits a non-trivial
$\C_+$-action then $(X,p)$ is rational. In particular, if
$X$ admits an  effective action of a non-abelian algebraic
group then $(X,p)$ is rational.
\ecor

\proof
The first statement is an immediate consequence of the
preceding theorem since the general orbits of the
$C_+$-action are closed rational curves in $X$ not passing
through $p$.

If $G$ is a non-abelian algebraic group
acting on $X$ then it contains an algebraic subgroup isomorphic to
$\C_+$. Indeed, if the solvable radical of
$G$ does not contain a subgroup isomorphic to $\C_+$ then
this radical is isomorphic to $\C^{*n}$ for some $n\ge 0$
and so $G$ is reductive. It follows then from the
structure theory of reductive groups that $G$ contains an algebraic
subgroup isomorphic to $\C_+$. Restricing the action of
$G$ to such a subgroup the second statement of the
corollary follows from the first one.
\qed

\brem
Note that \ref{newthm0} is no longer true if the
singularity is non-isolated as is seen by the
example of a product $X\times\pP^1$, where $X$ is a normal
variety with a non-rational point.
\erem

\subsection{Quasirational singularities of
quasihomogeneous varieties}

\bnota\label{qrv}
Let $A=\bigoplus_{\nu\ge 0} A_{\nu}$ be a
graded algebra of finite type over $\C$, and
let $V=\Spec A$ be the corresponding affine variety
with a good $\C^*$-action. We assume in this subsection that
the vertex set $V(A_+)$ (where $A_+:=\bigoplus_{\nu>0}
A_{\nu}$) consists of one point $p$
which is an isolated normal singularity.
\enota

Following Abhyankar \cite{Ab} in the 2-dimensional case we
introduce the following notion.

\bdefi\label{qr}
Consider the quotient variety
$\G:=(V\backslash \{p\})/\C^*=\Proj A$. The
singularity $(V,\,p)$
will be called {\it quasirational} if
$\G$ is uniruled\footnote{i.e., there is a rational dominant map
from $\pP^1\times Y\rmap X$ for some variety $Y$ or,
equivalently, for a general point $x\in X$ there is a rational
curve $C\subseteq X$ through $x$; see \cite[Chapt.\ IV]{Kol}.}.
\edefi

It is well known (see e.g.\ \ref{forms}) that in this case the
Kodaira dimension of $\G$ is $-\infty$. It is an open question
whether  conversely $k(\G)=-\infty$ implies that $\G$ is
uniruled; cf.\
\cite[IV 1.12]{Kol}.

Note that $\G$ has at most
cyclic quotient singularities
(see \cite{Fl} or \cite[Prop.\ 2.8]{FlZa1}). Hence
$H^0(\G,\cL_{\G}^{m})\cong H^0(\G,\omega_{\G}^{[ m]})$, where
$\omega_{\G}^{[m]}$ denotes the reflexive hull of
$\omega_{\G}^{\otimes m}$.

If $V$ is a surface then the condition of \ref{qr} means
that $\G$ is a rational curve or,
equivalently, that every irreducible
component of the exceptional divisor $\sigma^{-1}(p)$ of
the minimal resolution of singularity $\sigma:V'\to V$
is a rational curve (see e.g., \cite{Ab} and also
\cite{Or}; alternatively this follows from the fact
that the weighted blowing up of $V$ has at most quotient
singularities, see \cite{Fl}).

The following simple lemma provides a restriction
on a singularity
to be quasirational.

\blem\label{quasirat-forms}
If $(V,p)$ is quasirational then $(\omega_A)_0=0$.
In dimension 2 the
converse is also true.
\elem

\proof
By \cite[8.9]{Fl2} there is an exact sequence
$$
0\to (\Omega^1_\G)^{\vee\vee}\to
\tilde\Omega^1_A \stackrel{\xi}{\lto}
\cO_\G\to 0,
$$
where $\xi$ is induced by the Euler derivation and
$\tilde{\phantom{\G}}$ denotes the associated sheaf. Taking
determinants it follows that $\omega_\G\cong \tilde\omega_A$
and so $(\omega_A)_0\cong H^0(\G,\omega_{\G})$. As mentioned
above $\G$ has at most quotient singularities. Hence, if the
Kodaira dimension
$k(\G)=-\infty$ then
$H^0(\G,\omega_{\G})$ and
therefore also $(\omega_A)_0$ vanishes.
Conversely, if $\dim V=2$ then
$\G$ is a smooth complete curve, and so
the vanishing of
$H^0(\G,\omega_{\G})$ implies that this
curve is rational.
\qed

By an {\it orbit closure} we mean the closure in
$V$ of a one-dimensional
orbit of the $\C^*$-action; this is a rational curve
passing
through the vertex $p$ with a  normalization isomorphic
to the affine line $\C$.
We are interested in the existence of closed
rational curves
in $V$ different from orbit closures.
We have the following proposition.

\bprop\label{quasi} For a quasihomogeneous variety $V$ as in \ref{qrv}
the following
conditions are equivalent:

\begin{enumerate}
\item[(i)] $(V,\,p)$ is a quasirational singularity;
\item[(ii)] for a general point $x\in V$ there exists a closed rational
curve $C\subseteq V$ which
passes through $x$, with $C$
being different
from an orbit closure.
\end{enumerate}

\no If $\dim V=2$ these are also equivalent to:

\begin{enumerate}
\item[(iii)] The surface $V$ is rational.
\end{enumerate}
\eprop

\proof
Over an open dense subset, say,
$U\subseteq \G$ the canonical map
$\pi:V\backslash \{p\}\to\G$ is a trivial
$\C^*$-fibration so that
$\pi^{-1}(U)\cong U\times \C^*$. Assume first
that  (ii) is satisfied.
Projecting curves as in (ii) we see that for a general point $x$ of
$\G$ there is a rational curve in $\G$ passing through $x$,
whence $\G$ is uniruled.

Conversely, if $\G$ is uniruled then for a general point $x$ of $V$ there is a
rational curve
$C$ in $\G$ through $\pi(x)\in\G$. For any point $t\in\C^*$
then $(C\cap U)\times\{t\}$ is a rational curve in $U\times
\C^*\cong \pi^{-1}(U)$. Thus
$\pi^{-1}(U)$ is covered by
rational curves that are not contained in orbit closures
and (ii) follows.  For the equivalence
of (iii) with (i), (ii) in the two dimensional
case we refer the reader
to \cite[L. 7]{KalZa}.
\qed

  \brem\label{burns}
(1) According to \cite{Bur}
the affine cone $V$ over a smooth
cubic 3-fold $\G\subseteq \pP^4$ has a rational
and quasirational singularity at the origin, whereas
$\G=V^*/\C^*$ is not rational \cite{CleGri}.
Hence  in higher dimensions
the equivalence (i)$\Leftrightarrow$(iii)
does not hold, in general.

(2)
Clearly, a quasirational singularity of a variety
with a good $\C^*$-action is not necessarily rational.
For instance, the Pham-Brieskorn surface
$$
V_{p,q,r}\subseteq \C^3\,:\qquad
x^{p}+y^{q}+z^{r}=0
$$
has a rational singularity at the
origin only for the
Platonic surfaces
$V_{2,2,m}$, $V_{2,3,3}$, $V_{2,3,4}$, $V_{2,3,5}$ (see
\ref{cor ci}), whereas it is
quasirational in many other cases (namely, iff one of the
conditions (i), (ii) of
\ref{PHBR}(a) is fulfilled). The simplest example of a
quasirational non-rational singularity is
$(V_{2,3,7},\,{\bar 0})$
(see also \cite{Ab, Or} for further examples).
\erem

\bexa  Consider the surface
$V=\{F_d(x,y)-z^m\}\subseteq \C^3$, where $F_d$ is a
homogeneous polynomial of degree $d$ without multiple
factors. The singularity at the origin
of the surface $V$ is quasirational if and only if either
$d=2$ or $\gcd (m,d)=1$ \cite{KalZa}. At the same
time,  this singularity is rational
if and only if  $d\le 2$ or $(d,m)=(3,2)$, see Corollary
\ref{schm} below.
\eexa

\subsection{Polynomial curves on
quasihomogeneous varieties}

Recall that by a {\em polynomial curve} on a
variety $V$ we mean a rational curve
$C\subseteq V$ with one place at infinity,
i.e., a curve
$C$ with
a normalization
$C'\cong\C$.

If $V\subseteq \C^n$ is an affine variety
then
a normalization morphism
$\nu:\C\to C\hookrightarrow
V $ is given (in the coordinates of
the ambient affine space $\C^n$) by
a sequence of
$n$ polynomials of one variable that are not all constant;
this explains our terminology.

As a consequence of
\ref{dell} (b), (c), \ref{forms}(b) and \ref{Kempf} we obtain
the following statement.

\bprop\label{logkod} Let $V$ be a quasihomogeneous
normal algebraic variety \footnote{That is,  $V=\Spec A$ where
$A=\bigoplus_{i\ge 0}A_i$ is a normal
graded $\C$-algebra of finite
type.} with at most isolated
singularities. Suppose that for a point $p\in \mbox{Sing}\, V$ and
a general point
$x$ of $V$  there exists a polynomial curve
$C\subseteq V\backslash \{p\}$ through $x$. Then
$\delta_m(V,p)=0$
for all $m\ge 1$. In particular, $(V,p)$
is a rational singularity provided that it
is Cohen-Macaulay.
\eprop

The next corollary gives a
geometric application to $\C_+$-actions
on quasihomogeneous varieties (cf. Corollary \ref{agp}).

\bcor\label{cstar-action} Let $V$ be a quasihomogeneous
normal  variety with at most isolated
singularities and $p\in\Sing V$. If $V$ admits a non-trivial
$\C_+$-action then
the conclusions of \ref{logkod} hold. In particular, this holds
for any normal affine algebraic variety $V$ with isolated
singularities that admits an effective action of a non-abelian
non-nilpotent group
$G$.
\ecor

\proof This is an immediate consequence of \ref{logkod}.
Indeed, in the first case the orbit $\C_+.x$ through a
general point $x\in V$ is a polynomial curve in
$V\backslash
\{p\}$. The second statement follows from the first one since
$G$ being non-nilpotent contains an algebraic subgroup
isomorphic to $ \C^*$, whence $V$ is quasihomogeneous. As $G$ is
non-abelian it also has an algebraic subgroup isomorphic to
$\C_+$.
\qed

\bsit A similar result holds if $V$ admits a second
$\C^*$-action.
Of course, one can
obtain a second grading of $A$ (and hence
a second $\C^*$-action on $V$)
by simply multiplying degrees with a
constant $k$, i.e.\ $A=\bigoplus_{i\ge 0}A'_i$, where
$A_{kj}'=A_{j}$ and $A_{i}'=0$ for $i\not\equiv 0\mod k$.
We say that two gradings (resp., two $\C^*$-actions)
are {\em truly
different} if they cannot be obtained
by multiplying degrees
of a third one with appropriate constants.
Equivalently, this means that the associated
derivations span a $\C$-vector
space of dimension 2  in $\text{Der} A$.
\esit

\bthm\label{cplus-action}
Let $A=\bigoplus_{i\ge 0}A_i$ be a normal graded
$\C$-algebra of finite
type defining a quasihomogeneous variety $V=\Spec A$ with
at most isolated
singularities. If $V$ admits a second truly
different $\C^*$-action
then $\bar k(V\backslash  \Sing V)=-\infty$ and
$\delta_m(V,p)=0$ for all $m\ge 1$ and $ p\in\Sing V$.
In particular, $(V,p)$ is a log-terminal singularity provided
that it is $\Q$-Gorenstein.
\ethm

\proof
If $A_0\ne \C$ then the result follows from
\ref{dell} (b) and \cite[1.26, 2.27]{FlZa1}.
Suppose now that $A_0=\C$.  Let $\delta$ be the
derivation associated to the second $\C^*$-action
and write $\delta=\sum_{i=k}^l \delta_i$,
where $\delta_i$
is a homogeneous derivation of
degree $i$. We may assume that $\delta_k$,
$\delta_l\ne 0$. If $k<0$ then
$\delta_k$ is locally nilpotent and so defines
(see \cite{Ren}) a
$\C_+$-action. Applying
\ref{cstar-action} we are done. Thus we may
suppose that $k\ge 0$.

Choose a finite homogeneous generating set, say, $E$
for $A$ as a $\C$-algebra. As $\C^*$ is
reductive there is a
$\C^*$-invariant subspace
$V\subseteq A$ containing $E$ with $\dim_{\C} V <
\infty$. For $N$ sufficiently large $V$
is contained in
$\bigoplus_{i=0}^N A_i$, and so
$\delta^n(a)\in \bigoplus_{i=0}^N
A_i$ for every $a\in E$ and $n\in \N$. If $l>0$
then (considering the leading term of
$\delta^n(a)$) it follows that
$\delta_l^n(a)=0$ for $n\ge N$. Hence $\delta_l$ is
locally nilpotent, and we can
again conclude by \ref{cstar-action}.

Let us finally treat the case where $\delta=\delta_0$
is homogeneous of
degree 0. In this case the two $\C^*$-actions
commute, and so they
define a bigrading on
$A$. In particular there is a grading of $A$ with $A_0\ne \C$.
Applying the first case the desired result follows.
\qed

\section{Rational curves on quasihomogeneous surfaces}

\subsection{Polynomial curves and quotient singularities}

In the surface case,
Proposition \ref{logkod} can be strengthened as follows.

\bthm\label{newthm} Let $V$ be an affine
normal algebraic surface
and let $p\in V$ be a singular point.
Assume that for a general point $x\in V$ there exists a
polynomial curve
$C\subseteq  V$ with $x\in C$ but not  passing through $p$.
Then $(V,p)$ is a  quotient singularity.
\ethm

We need  the following lemma.

\blem \label{newlem1}
If a normal surface singularity $(V,p)$ is
finitely dominated by a quotient singularity $(W,q)$
then it is also a quotient singularity.
\elem

\proof
By \cite[Cor.\ 1.27]{FlZa1} we have that $\delta_m(V,p)\le
\delta_m(W,q)=0$ for all $m\ge 1$. Thus by
\ref{Watanabe}
$(V,p)$ is a quotient singularity.
\qed

The proof of Theorem \ref{newthm} is
based on Proposition \ref{newlem2} (b) below. Recall that
a complex surface $\Sigma$ is said to be {\it ruled}
(resp., {\it affine ruled})
if it is equipped with a proper morphism (resp., a morphism)
  $\pi :\Sigma\to S$ (called {\it ruling})
onto a smooth curve $S$
with general fibers isomorphic to $\pP^1$ (resp., to $\C$).

\bprop\label{newlem2}
Let $\pi:\Sigma\to S$ be a normal ruled surface over a smooth curve $S$
with
generic fiber $F\cong\pP^1$. Then the following hold. If there is
a Weil divisor
$H$ on $\Sigma$ with $HF=1$ such that $\Sigma\backslash H$
is affine, then $\Sigma\backslash H$ has at most cyclic
quotient singularities.
\eprop

This  is an immediate consequence of the following result due
to Miyanishi  \cite{Miy} (see also \cite{Gur}).

\bthm {\rm (Miyanishi)} \label{Miy}
Let $X$ be a normal affine algebraic surface, and assume that
$X$  contains a cylinder, i.e. a Zariski
open dense subset $U\cong \C\times S_0$
where $S_0$ is a curve. Then $X$ has
at most cyclic quotient singularities.
\ethm

\no {\it Proof of Theorem \ref{newthm}.}
Let $\bar V$ be a projective compactification of $V$ and denote
$D:=\bar V\backslash V$ the divisor at infinity.
As in the proof of Proposition \ref{forms},
there exists a projective normal ruled surface $\pi:\Sigma\to S$ over a
smooth complete curve $S$ and a surjective morphism $f: \Sigma\to
\bar V$  such that
$\pi\times f:\Sigma\to S\times \bar V$ is finite over a closed subscheme
of $S\times \bar V$ and such that the divisor $E:=f^{-1}(D)_\red$
on $\Sigma$ intersects the general fiber of
$\pi$ in just one point.  Consider now $p\in V$ and let
$q\in\Sigma$ be a point over $p$. By construction $q$ is not
contained in $E$. For a fiber $F_s$ over a point $s$ of $S$ different
from $\pi(q)$ the scheme
$\Sigma\backslash (E\cup F)$ is finite over a closed subscheme of the
affine scheme
$V\times (S\backslash\{s\})$ and so it is affine as well. Applying
\ref{newlem2}(b) to the divisor $H:=E\cup F_s$ we obtain that
$(\Sigma,q)$ is at most a cyclic quotient singularity. As
$(V,p)$ is dominated by $(\Sigma, q)$ under a (analytically) finite map,
the assertion follows from Lemma \ref{newlem1}.
\qed

Note that the assumption that $V$ is affine in Theorem
\ref{newthm} is essential. For instance, blowing up a
ruled surface repeatedly in one fibre and blowing down any
connected curve in the exceptional locus results in a surface
with an isolated singularity and many polynomial curves.
However,  the minimal resolution of a quotient singularity is
star shaped (in fact, this is true for any isolated quasihomogeneous surface
singularity, see \cite{OW}) whereas by the procedure above we can
easily produce singular points, whose minimal resolution is not
star shaped.

In view of \ref{newthm} and \ref{Watanabe} it is natural to ask
the following question.

\begin{que}
Let $V$ be an affine
normal algebraic variety
and let $p\in V$ be an isolated singular point.
Assume that for a general point $x\in V$ there exists a
polynomial curve
$C\subseteq  V$ with $x\in C$ but not  passing through $p$.
Is then $\delta_m(V,p)=0$ for all $m\ge 1$?
\end{que}

\subsection{Polynomial curves on quasihomogeneous surfaces}
Recall that a {\it small finite
subgroup}
$G\subseteq {\rm GL}\,(2,\,\C)$ is
a finite subgroup which does not contain any
pseudo-reflection, that is,
a linear transformation of finite order $\sigma\in {\rm
GL}\,(n,\,\C)\backslash\{{\bf 1}_n\}$  which is the identity
on a hyperplane.
The following fact is well known; in lack of a reference we
include a brief argument.

\blem\label{glob} A surface $V$ with a
good $\C^*$-action
has an isolated quotient singularity at the vertex
$p\in V$ if and only if $V\cong\C^2/G$
where $G\subseteq {\rm GL}\,(2,\,\C)$ is a small finite
subgroup.\elem

\proof
  From the assumption that $(V,\,p)$ is
an isolated quotient singularity
it follows that the fundamental group
$G:=\pi_1\,(V\backslash \{p\})$
(isomorphic to the local fundamental group of the
singularity $(V,\,p)$)
is finite. Thus
the universal cover $W'\to V\backslash \{p\}$
extends to a branched covering
$\tau:W\to V$ such that $W$ is a smooth surface
and
$\tau^{-1}(p)$ consists of a single point, say
$q\in W$.
After replacing the original $\C^*$-action on $V$
by the new one\footnote{Hereafter we denote by $t.v$
the image of
a point $v\in V$ under an element $t$ of a group acting on $V$.}
$t^{\ord G}. v$, where $t\in\C^*$, $v\in V$, it can be lifted to a good
$\C^*$-action on $W$
with a unique fixed point $q\in W$,
so that the morphism $\tau$ becomes equivariant.
It is well known that in this situation $W\cong\C^2$
and, what is more, $B:=H^0(W,\,{\mathcal O}_W)=\C[X,Y]$ where
$X,\,Y\in B$ are two algebraically independent
homogeneous elements of minimal possible degrees, say,
$d_X$ resp., $d_Y$ (geometrically, this isomorphism is
composed of an embedding
$W\hookrightarrow \C^N$ equivariant with respect
to a diagonal linear $\C^*$-action on
$\C^N$, followed by an equivariant projection onto
the tangent plane
at the fixed point $q\in W$;
see \cite[8.5]{Za}). Thus also the quotient morphism
$\tau:W\to V$ is equivariant with respect to the
$\C^*$-action $t^d. v$
on $V$ and the diagonal linear $\C^*$-action
$$\lambda: (X,\,Y)\longmapsto
(\lambda^{d_X}\cdot X,\,\lambda^{d_Y}\cdot Y) $$
on $\C^2=\Spec \C[X,Y]$.
This completes the proof.
\qed

If $G\subseteq {\rm GL}\,(2,\,\C)$ is a small
finite subgroup and $\tau:\C^2\to V:= \C^2/G$
is the quotient map
then there exists a polynomial curve
$\C\to V\backslash\{p\}$,
where $p=\tau({0})\in V$.
In the next theorem we show that the converse is also true.

\bthm\label{mainthm} Let $V$ be a normal affine surface
with a good $\C^*$-action and an isolated singularity
$p\in V$. Then the following conditions are equivalent:

\begin{enumerate}
\item[(i)] ${\overline k}\,(V\backslash \{p\})=-\infty$;

\item[(ii)] $\delta_m\,(V,\,p)= 0$ for all $m\ge 1$;

\item[(iii)] $(V,\,p)$ is a quotient singularity;

\item[(iv)] there exists a polynomial curve
$\C\to V\backslash\{p\}$.
\end{enumerate}
\ethm

\proof
The equivalences (i)$\Leftrightarrow$(ii)
and (ii)$\Leftrightarrow$(iii)
follow from Theorem \ref{dell} (c), resp.\ Watanabe's theorem
\ref{Watanabe}. The implication (iii)$\Rightarrow$(iv)
has been observed above. For the proof of (iv)$\Rightarrow$(i)
note that moving a polynomial curve in $V\backslash \{p\}$
with the $\C^*$-action provides polynomial curves through a
general point of $V\backslash \{p\}$. Thus this implication
follows from Proposition
\ref{forms} (b).
\qed

\brem
The implication (i)$\Rightarrow$(iii) provides a strengthening
of Proposition 4.7
in\footnote{Cf. also \cite{MaMiy}.} \cite{Be1}
which states that a
quasihomogeneous normal surface singularity
$(V,\,p)$ with ${\overline
k}\,(V\backslash \{p\})=-\infty$ is rational.
Recall that any quotient
singularity is rational (cf.\ \cite{Bri, Bur, Vie}).
Moreover, the quotient $V/G$ of an isolated
rational
singularity $(V,\,p)$ by a finite group $G$
acting freely off the
singular point $p\in V$  is again rational
  (cf. Lemma \ref{newlem1}).
\erem

On the other hand,
it is well known that a rational singularity
$(V,\,p)$ of  a surface $V$  with a good $\C^*$-action
is not necessarily a quotient singularity
(and if so then by \cite[4.7]{Be1} or
by Theorem \ref{mainthm}
${\overline k}\,(V\backslash \{p\})\ge 0$).
To construct concrete examples we use the following
proposition.

\bprop\la{quot exa}
Let $V$ be an affine normal surface with a good
$\C^*$-action and a quasirational singularity at its vertex,
say $p\in V$. Let $\Z_d\cong \langle\zeta\rangle$ be the
subgroup of $\C^*$ generated by a primitive $d$-th root of unity
$\zeta$. Then the following hold.

1. If $d\gg 0$ then the quotient
$W:=V/\Z_d$ has at most a rational singularity.

2. ${\bar k}\,(W\backslash \{q\})={\bar k}\,
(V\backslash \{p\})$, where $q$ denotes the vertex of $W$. In
particular,
$W$ has a quotient singularity if and only if $V$ has a quotient
singularity.
\eprop

\proof
Let $A$ be the coordinate ring of $V$ with its natural
grading $A=\bigoplus_{i\ge 0}A_i$ induced by the
$\C^*$-action. The coordinate ring of $W$ is then given
by the Veronese subring $B:=A^{\Z_d}=\bigoplus_{i\ge 0} A_{id}$.
Moreover
$$
\omega_B^{[m]}\cong
(\omega_A^{[m]})^{\Z_d}=\bigoplus_{i\in\Z}(\omega_A^{[m]})_{id}
$$
for all $m\ge 1$. Now  $(\omega_A)_{\le 0}$
is a vector space over $\C$ of finite dimension and
$(\omega_A)_0=0$, since $V$ has a quasirational singularity (see
\ref{quasirat-forms}). Hence for $d\gg 0$ the homogeneous components
$(\omega_A)_{id}$ vanish for all $i\le 0$ and so
$\omega_B=(\omega_B)_{>0}$. Using \ref{dell}(a) it follows that
$W$ has a rational singularity at its vertex $q\in W$, proving
(1).

To deduce (2), note that by \ref{dell}(c)
$$
\bar p_m(V\backslash\{p\})=\dim_\C (\omega_A^{[m]})_0
\quad\mbox{and}\quad
\bar p_m(W\backslash\{q\})=\dim_\C (\omega_B^{[m]})_0.
$$
As $(\omega_B^{[m]})_0\cong
(\omega_A^{[m]})_0$ for $m\ge 1$ this implies that
$V\backslash\{p\}$ and $W\backslash\{q\}$ have the same
logarithmic Kodaira dimension. The second part of (2) is now a
consequence of \ref{mainthm}.
\qed

\bexa\label{nonquo}
As a concrete example, consider the Pham-Brieskorn
surface $V=V_{p,q,r}$ given by the equation
$x^p+y^q+z^r=0$ in $\C^3$. This surface admits a $\C^*$-action via
$$
t.(x,y,z):=(t^{qr}x,t^{pr}y,t^{pq}z).
$$
Assume that $1/p+1/q+1/r<1$  and that
one of the conditions (i), (ii) of Theorem \ref{PHBR} (a)
is fulfilled so that $V$ has a
non-rational quasirational singularity at the origin. If
$\zeta $ is a $d$-th root of unity and $d\gg 0$ then by
Proposition \ref{quot exa} the
quotient $W:=V/\langle \zeta\rangle $ has a rational singularity
which is not a quotient one. Using Corollary \ref{cor ci} the
logarithmic Kodaira dimension of
$V\backslash\{p\}$ (and hence also of $W\backslash\{q\}$) is
equal to 1.

Notice that by Theorem \ref{mainthm}, both $W\backslash \{q\}$ and
$V\backslash
\{p\}$ do not contain polynomial curves,
while by Proposition \ref{quasi} $V$ and $W$ contain closed
rational curves different from an orbit closure. Clearly
$V\setminus \{p\}$ cannot contain closed rational curves, see
\ref{newthm0}. However, we do not know whether
$W\setminus
\{q\}$ contains closed rational curves. More generally one is
tempted to pose the following question.
\eexa

\begin{que}
Does an affine normal surface $V$ with a good $\C^*$-action and
a rational singularity contain a closed
rational curve not passing through its vertex?
\end{que}

\subsection{Surfaces with $\C_+$-actions and cyclic
quotient singularities}
Let $V$ be an affine normal surface with a good $\C^*$-action.
If $V$ admits a $\C_+$-action then by Theorem \ref{mainthm}
it has at most a quotient singularity. In this
subsection we examine the question as to when a
quasihomogeneous surface with a quotient singularity admits a
$\C_+$-action.

\bsit
Recall that a
{\em cyclic quotient surface singularity} arises
locally as the quotient of $\C^2$ by the action
of a group $\Z_d\cong \langle \zeta\rangle$ of $d$-th roots of unity
generated by
\be\label{act} (x,\,y)\longmapsto
(\zeta x,\,\zeta^e y),\ee where
$\zeta$ is a primitive
$d$-th root of unity and $e\in \N$, $\gcd (e,d)=1$.

By a theorem of Brieskorn \cite{Bri},
a normal surface singularity
$(W,\,p)$ is a cyclic quotient singularity if and only if
the dual graph $\G$ of
its minimal embedded resolution $\sigma:W'\to W$ is linear
(this linear weighted graph $\G$ represents a
{\it Hirzebruch-Jung string}; see e.g.,
\cite[III.2, III.5]{La1, BPVdV}).
By \cite{KKMSD}
(see also \cite{Pi1}), $(W,\,p)$ is a cyclic
quotient singularity
if and only if it is a toric surface
singularity, that is, $W$
is an affine surface with a good effective
$\C^*\times \C^*$-action
and $p$ is its only fixed point.
\esit

We characterize below affine
surfaces with
a good $\C^*$-action
which  also admit a regular action
of the additive group $\C_+$
(cf.\ also \cite{DanGi, Be2}).

\bthm\label{cycl} Let $V$ be an affine surface
with a good $\C^*$-action and an isolated normal
singularity $p\in V$. The surface $V$ admits a non-trivial
$\C_+$-action if an only if $(V,\,p)$ is
a cyclic quotient singularity:
$V\cong\C^2/G$ with $G\cong \Z/d\Z$
being a small finite subgroup of
${\rm GL}\,(2;\,\C)$. \ethm

The proof is given in Lemmas
\ref{cyc}--\ref{cyquo}
below.

\blem\label{cyc} Let $V=\C^2/G$
be a quotient of $\C^2$
with respect to a small subgroup
$G\subseteq{\rm GL}\,(2;\,\C)$.
If $G$ is a cyclic group
then the surface $V$ admits a nontrivial
$\C_+$-action. \elem

\proof We may assume that the action of the group
$G\cong\Z_d$ on $\C^2$ is generated by the linear transformation
\begin{displaymath}
\mbox{$g=$}
\left( \begin{array}{cc}
\zeta & 0\\
0 & \zeta^e
\end{array} \right)\,,
\end{displaymath}
where $\zeta$ is a primitive $d$-th
root of unity. As $G$ is small,
$e$ and $d$ have to be coprime.
The locally nilpotent derivation
$$\partial:\C[X,Y]\longrightarrow \C[X,Y]\,,$$
$$\partial(X):=0,\qquad \partial(Y):=X^e$$
commutes with the action of $G$ i.e.,
$\partial(g. f) = g. \partial(f)$
for every $f\in \C[X,Y]$ and every $g\in G$.
Hence $\partial$
descends to a non-zero
locally nilpotent derivation say, $\partial'$
on the subalgebra of invariants $A:= \C[X,Y]^G$
and so $\C_+$ acts effectively on $A$
(as well as on the surface $V=\Spec A$)
by $t\longmapsto \exp\,(t\partial')$
(see e.g., \cite{Ren}).\qed

The converse is also true (see Lemma \ref{cyquo} below).
But first we need the next simple lemma. Recall that
a derivation, say, $\p$ of a domain, say, $A$
(over $\C$) naturally extends to the fraction field
$\Frac A$ leaving invariant the integral
closure $\bar A$ of $A$ in $\Frac A$
\cite{Sei}.
We have the following statements.

\blem\label{sile} Let $A\hookrightarrow B$
be an integral
extension of finitely generated domains over $\C$,
and let $\partial\in \Der\, B$
be a derivation
such that $\partial(A)\subseteq A$. Then the following hold.

\begin{enumerate}
\item[(a)] \cite{Vas} If the restriction
$\partial\,|\,A$ is locally nilpotent  then so is $\partial$.
\item[(b)] Suppose moreover that $B$ is graded and that $A$ is a
graded subring.  If the restriction
$\partial\,|\,A$ is homogeneous then $\partial$ is
homogeneous.
\end{enumerate}\elem

\proof (a) We provide an alternative geometric
proof supposing that
$A$ and $B$ are defined over $\C$. Let $\pi:Y:=\Spec B\to X:=\Spec
A$ be the map induced by $A\hto B$. Every orbit
$\orbit_x:=\C_+.x\subseteq X$, where $x\in X$ is a closed point, of
the associated $\C_+$-action on $X$ is isomorphic to $\C_+$ or reduced
to a point. As
$X$ is affine, $\orbit_x$ is closed in $X$. By assumption the
infinitesimal generator of the $\C_+$-action lifts to a vector
field on $Y$ and so $\pi^{-1}(\orbit_x)\to\orbit_x$ is unramified. Thus
the preimage
$\pi^{-1}(\orbit_x)$ is a disjoint union of copies of $\C_+$ or of
points, and the
$\C_+$-action on $\orbit_x$ lifts to an (algebraic) action on
$\pi^{-1}(\orbit_x)$. In particular, the derivation $\partial$ induces
a locally nilpotent derivation on the affine coordinate ring of
$\pi^{-1}(\orbit_x)$.

For an algebraic function $f\in B$ consider the subset
$A_n$ of all closed points $x\in X$ such
that the restriction of $\partial^n(f)$ to $\pi^{-1}(\orbit_x)$
vanishes identically. By the above reasoning, $\bigcup_{n\ge 0}A_n$
is the set of all closed points of $X$. Thus for some $n$ the
Zariski closure of $A_n$ in $X$ must be equal to $X$, and then
$\partial^n(f)$ vanishes identically on $Y$.

(b) Let $f\in B_k$ be a homogeneous
element and consider
an equation of integral dependence
$P(f)=0$, where
$$
P=X^n+a_{1}X^{n-1}+\dots +a_n
$$
with $a_i\in A$. We may assume that
$a_{\nu}$ is homogeneous of degree $k\nu$.
Applying $\partial$ gives
$$0=\partial(P(f))=P'(f)\cdot
\partial(f)+P_{\partial}(f)
\qquad
{\rm with}\qquad
P_{\partial}:=\sum_{\nu=0}^n
{\partial}(a_{n-\nu})X^{\nu}\,.$$
As $P_{\partial}(f)$ and $P'(f)$
are both homogeneous and $B$
is a domain, it follows that
${\partial}(f)$ is also homogeneous. \qed

\blem\label{cyquo} Let  $V$ and $p\in V$ be
as in \ref{cycl}. If  $V$ admits a non-trivial
$\C_+$-action then $V\cong \C^2/G$  for a small cyclic subgroup
$G\subseteq {\rm GL}\,(2;\,\C)$. \elem

\proof Since $V$ admits a non-trivial $\C_+$-action
with a fixed point $p\in V$, we have
${\overline k}\,(V\backslash \{p\})=-\infty$,
and hence
by Theorem \ref{mainthm}, $(V,\,p)$
is a quotient singularity.
Thus,
Lemma \ref{glob} can be applied.
We keep the notation as
in the proof of this lemma;
in particular, we have $V=W/G$,
where $W\cong \C^2$ and
$G\subseteq {\rm GL}\,(2,\,\C)$ is
a small finite subgroup, with
quotient morphism
$\tau:W\to V$. Furthermore,
$B:=H^0(W,\,{\mathcal O}_W)\cong\C[X,Y]$ with
$X,\,Y\in B$ being homogeneous of
degree $d_X$ resp., $d_Y$,
and $\tau$ is equivariant with
respect to the corresponding
(linear diagonal) $\C^*$-action
on $\C^2=\Spec \C[X,Y]$.

The derivation $\partial=\partial_{\varphi}\in
\Der\, A$ corresponding to the
$\C_+$-action, say $\varphi:\C_+\times V \to V$ is known to be locally
nilpotent (see e.g., \cite{Ren}) and non zero.
Since the leading term $\partial_\ell$ of
$\partial$ is again a locally
nilpotent derivation (see e.g., \cite{Ren, ML})
we can replace
$\partial$ by $\partial_\ell$
and the action
$\varphi(t,-)=e^{t\partial}$ by the
new action $e^{t\partial_\ell}$.
Thus we may assume that
$\partial$ is homogeneous.

The derivation $\partial$ is induced by a homomorphism, say,
$h:\Omega^1_A\to A$ and so it gives a homomorphism $h\otimes
1:\Omega^1_A\otimes_AB\to B$. As $A\to B$ is unramified in codimension
1, the reflexive hull of  $\Omega^1_A\otimes_AB$  is just $\Omega^1_B$,
and so we get an induced homomorphism $\tilde h:\Omega^1_B\to B$.
Composing this map with the total differential $d:B\to \Omega^1_B$ we
obtain a derivation $\tilde\partial$ that lifts $\partial$ and is
locally nilpotent by \ref{sile}(a). Taking the exponential of
$\bar\partial$ the
$\C_+$-action
$\C_+\times V \to V$ lifts to an action
$\C_+\times W \to W$ such that the diagram
\begin{diagram}
\C_+\times W & \rTo^{\tilde \varphi} & W\\
\dTo<{\id\times\tau} && \dTo>\tau\\
\C_+\times V & \rTo^{\varphi} & V
\end{diagram}
is commutative. Let us show that this action commutes with the action
of $G$ on $W$.
Indeed, it preserves the $G$-orbits and so, given $z\in\C^2$ and
$g\in G$, for any $\lambda \in\C_+$
there exists an element ${\tilde g}(\lambda)\in G$ such that
$\lambda . (g. z)={\tilde g}(\lambda).(\lambda . z)$.
Since ${\tilde g}(\lambda)$ is a continuous function on
the connected variety $\C_+$
with values in $G$, it must be constant:
${\tilde g}(\lambda)={\tilde g}(0)=g$, and so
$g\lambda=\lambda g$ for all
$g\in G,\,\,\lambda\in\C_+$,
as stated.

In what follows we distinguish two cases:

\smallskip\no {\it Case} (1): $d_X<d_Y$.
As the action of $G$ on the algebra
$B\cong\C[X,Y]$ is homogeneous, it follows that
$$g. X=\alpha(g)\cdot X\qquad {\rm and}\qquad g.
Y=\beta(g)\cdot Y$$
for two characters $\alpha,\,\beta:G\to S^1$.
For
an element $g\in G$, let $\ord \alpha(g)=kp_1$
and $\ord \beta(g)=kp_2$ with $\gcd (p_1,p_2)=1$.
If $p_1\neq 1$ then $\alpha(g^{kp_1})=1$ and
$\beta(g^{kp_1})\neq 1$,
so that the element $g^{kp_1}$ acts on $W=\C^2$
as a pseudo-reflection.
This leads to a contradiction since
the covering $\tau:W\to V$
is unramified outside the origin
$q={\bar 0}\in\C^2=W$.

Hence it follows that $p_1=1$ and (by symmetry) $p_2=1$ so
that the orders of
$\alpha(g)$ and
$\beta(g)$ are equal
for any $g\in G$. Thus
$G\cong \alpha(G)\subseteq S^1$
is a cyclic group, as desired.

\smallskip\no {\it Case} (2): $d_X=d_Y$.
Let $\widetilde{\partial}=
\widetilde{\partial}_{\widetilde{\varphi}}$
be the infinitesimal generator of the $\C_+$-action
$\widetilde{\varphi}$ on $B$ so that
$\widetilde{\partial}$ is a locally nilpotent derivation  of
the polynomial algebra $B=\C[X,Y]$. By our assumption,
$\partial$ is homogeneous and so, by Lemma \ref{sile},
$\widetilde{\partial}$ is as well homogeneous
of degree, say $d$. As ${\widetilde{\partial}}$
is locally nilpotent, we can find a homogeneous polynomial
$P(X,Y)\in B$
of positive degree with ${\widetilde{\partial}}(P)=0$.
The kernel
$\ker{\widetilde{\partial}}$ of
${\widetilde{\partial}}$ being factorially closed
(see e.g., \cite{ML}),
every factor in the decomposition $P=\prod_{i=1}^m
(\alpha_iX+\beta_iY)$ vanishes under
${\widetilde{\partial}}$. Hence after a linear
change of coordinates
(which is homogeneous since $d_X=d_Y$) we may assume that
${\widetilde{\partial}}(X)=0$. Write now
$${\widetilde{\partial}}(Y)=\sum_{i=0}^s Y^iF_i(X)\in B_{d+d_Y}$$
with $F_s\not\equiv 0$.
If $s\ge 1$ then an easy induction shows that
$${\widetilde{\partial}}^n(Y)=c_n\cdot Y^{ns-n+1}F_s(X)^n+{\rm
(lower\,\,\,order\,\,\,terms\,\,\,in}\,\,\,Y{\rm )}\neq 0$$
with $c_n\in\N$, $c_n\neq 0$,
which contradicts the assumption that ${\widetilde{\partial}}$ is
locally nilpotent.  Hence $s=0$ and so
$${\widetilde{\partial}}(Y)=c\cdot X^{d+d_Y}$$ for some $c\in
\C^*$. Replacing $Y$ by $Y/c$ we can achieve that $c=1$.

As ${\widetilde{\partial}}$ commutes with the action of the group
$G$ on $\C^2$ we have ${\widetilde{\partial}}(g. X)=g.
{\widetilde{\partial}}(X)=0$ for every $g\in G$. It follows that
$g. X=\alpha(g)\cdot X$ for a character $\alpha:G\to S^1$.
Furthermore, if $g. Y=\xi X+\eta Y$ then the equalities
$$\alpha(g)^{d+d_Y}\cdot X^{d+d_Y}=g. X^{d+d_Y}=g.
{\widetilde{\partial}}(Y)=
{\widetilde{\partial}}(g. Y)=\eta
{\widetilde{\partial}}(Y)=\eta\cdot X^{d+d_Y}$$
show that $\eta=\alpha(g)^{d+d_Y}$. If $\alpha(g)=1$
then the matrix of $g$ has the form
\begin{displaymath}
\mbox{$g=$}
\left( \begin{array}{cc}
1 & \xi\\
0 & 1
\end{array} \right)\,.
\end{displaymath}
Since $g$ has a finite order, $\xi$ must be zero.
In other words, the homomorphism $\alpha:G\to S^1$ is injective
and so $G\cong \alpha(G)$ is a cyclic group.
This completes the proof of the lemma and of Theorem \ref{cycl}.
\qed

   The next example shows that the assumption
of normality in Theorem \ref{cycl} is essential. Indeed,
the existence of a $\C_+$-action on an isolated surface
singularity does not imply that this singularity is normal.

\bexa
Consider the polynomial algebra
$\C[X,\,Y]$ with the derivation $\partial:=X\cdot
\frac{\partial}{\partial Y}$. Apparently, $\partial$
is a locally nilpotent, homogeneous derivation of
degree $0$ with respect to the standard grading on
$\C[X,\,Y]$. Hence $\partial$ induces a homogeneous
locally nilpotent derivation on the subring
$$A:=\C[X^2,\,XY,\,Y^2,\,X^3,\,Y^3]=\bigoplus_{\nu\neq 1}
\C[X,\,Y]_{\nu}\,.$$
Note that the normalization of $A$ is $\C[X,\,Y]$ and that $A$
is not a quotient singularity.
\eexa

\subsection{Rational curves on surfaces with a
Gorenstein singularity}

Recall the following notion.

\bdefi
A singularity $(V,\,p)$ is called
{\it Gorenstein} if it is Cohen-Macaulay and
$\omega_{V,p}\cong\cO_{V,p}$.
\edefi

Note that a normal surface singularity is always Cohen-Macaulay, and
so it is Gorenstein if and only if there exists a nowhere vanishing
holomorphic 2-form $\omega$ on $V\backslash\{p\}$.

We recollect below some useful facts on
Gorenstein singularities.

\bthm\label{gor}
\begin{enumerate}
\item[(a)] {\em (See e.g.,\ \cite[A.2.5]{FlVo})}
Any complete intersection singularity is
Gorenstein.
\item[(b)] \cite{KeWa}
A quotient singularity $(V,p)\cong(\C^n/G,\,{\bar 0})$
with  a small
finite subgroup $G$ of ${\rm GL}\,(n;\,\C)$
is Gorenstein if and only if $G\subseteq {\rm SL}\,(n;\,\C)$.
\item[(c)] In particular, a
cyclic quotient surface singularity $(V,\,p)\cong (\C^2/G,\,{\bar 0})$
is Gorenstein if and only if $G\subseteq {\rm SL}\,(2,\,\C)$
is a cyclic subgroup generated by
$$(x,y)\longmapsto (\zeta x,\,\zeta^{-1} y)$$
where $\zeta\in \C$ is a primitive root of unity of degree,
say, $d$.
In the latter case, $(V,\,p)\cong (V_{2,2,d},\,{\bar 0})$
where $V_{2,2,d}\subseteq \C^3$
is the dihedral surface with the equation $x^2+y^2+z^d=0$.
\item[(d)] \cite[2.5, 4.6]{ Wah}
A rational surface singularity $(V,\,p)$
is Gorenstein if and only if it is a rational double point
(i.e.,
${\rm mult}_p\, V=2$). In the latter case,
$(V,\,p)$ is a quotient singularity
$(\C^2/G,\,0)$ with a small finite subgroup
$G\subseteq {\rm SL}\,(2;\,\C)$
(e.g., see \cite {Pi2}).
\item[(e)]  \cite[4.5, A.6]{Wah}
Every
rational surface singularity is $\Q$-Gorenstein, and so a cyclic quotient
(in a canonical way) of
a Gorenstein singularity (called {\em the canonical cover}).
The canonical cover
is rational if and only if the original
singularity is a quotient singularity.
\end{enumerate}\ethm

One can consult e.g., \cite{TSH}
and \ref{nonquo} above
for concrete examples
of affine surfaces with a good $\C^*$-action and
an isolated rational singularity which
is not Gorenstein.
Note that $(V_{p,q,r}, {\bar 0})$ in \ref{nonquo}
is just the canonical
Gorenstein cyclic cover of $(W,q)$.

For Gorenstein surface singularities
we have the following result.

\bthm\label{lines} Let $V$ be an affine surface
with a good $\C^*$-action and an isolated Gorenstein
singularity $p\in V$.
The following conditions are equivalent:

\begin{enumerate}
\item[(i)] ${\overline k}\,(V\backslash \{p\})=-\infty$;
\item[(ii)] $(V,\,p)$ is a quotient singularity;
\item[(iii)] $(V,\,p)$ is a rational singularity;
\item[(iv)] there exists a closed rational curve
$C\subseteq V\backslash \{p\}$;
\item[(v)] there exists a polynomial curve
$f:\C\to V\backslash \{p\}$.
\end{enumerate}\ethm

\proof The equivalences (i) $\Leftrightarrow$ (ii)
$\Leftrightarrow$ (v) have been established in
Theorem \ref{mainthm} above.
The implication (v) $\Rightarrow$ (iv) is evident;
(iv) $\Rightarrow$ (iii)
follows from Theorem \ref{newthm0}. For Gorenstein
singularities,
(ii) and (iii) are equivalent; see \cite{Bri}
for (ii) $\Rightarrow$ (iii)
and Theorem \ref{gor}(d) above for (iii) $\Rightarrow$ (ii).
This proves the theorem.
\qed

\bexa
The Pham-Brieskorn hypersurface
$$
V_{\bar p}:=\Big\lbrace\sum_{i=1}^n x_i^{p_i}=0\Big\rbrace
\subseteq \C^n\qquad (\mbox{with}\quad p_i\ge 2\quad\forall i)
$$
is quasihomogeneous with weights $w_i=1/
p_i$, $i=1,\dots,n$, and the defining equation has degree 1.
According to \ref{cor ci} $(V, 0)$
is a rational singularity if and only if
$$ \sum_{i=1}^n \frac{1}{p_i}>1\,.$$
On the other hand, it is known
\cite[Prop.5.2]{Ste, BrMa, DarGr}
that for
$$
  \sum_{i=1}^n \frac{1}{ p_i}\le \frac{1}{n-2}\,
$$
the Diophantine equation $\sum_{i=1}^n x_i^{p_i}=0$  does not
admit non-constant solutions in polynomials coprime in
pairs, that is  any polynomial curve $\C\to V_{\bar p}$
meets a codimension
2 coordinate subspace of $\C^n$.

In particular,  for $n=3$ (in virtue of
Theorem \ref{lines}) we have the equivalences
$$
{\overline k}\,(V_{\bar p}\backslash
\{\bar 0\})=-\infty\qquad \Leftrightarrow \qquad
\frac{1}{p_1}+\frac{1}{p_2}+\frac{1}{p_3}>1
$$
$\Leftrightarrow$
the surface $V_{p_1,p_2,p_3}$ has a rational
(or equivalently, a quotient) singularity
$\Leftrightarrow$
it is one of the Platonic surfaces
$V_{2,2,n},\,\,V_{2,3,3},\,\,V_{2,3,4},\,\,V_{2,3,5}$.
It has a cyclic quotient ${\rm (}A_{d-1}{\rm -)}$singularity
if and only if it is a dihedral
surface $V_{2,2,d}$
(cf. \cite[4.6]{Be1} and \cite{MaMiy}),
and only in this case it admits an effective $\C_+$-action
\cite[L. 4]{KalZa}.
\eexa

\bsit\label{qn} {\bf Question.} The affine Fermat cubic 3-fold
$$V:=\{x_1^3+\ldots+x_4^3=0\}$$ in $\C^4$ has a rational singularity at the
origin. Does it admit a non-trivial regular $\C_+$-action?\esit

Other examples are provided by the following result.

\bcor\label{schm} Let $V$ be a surface in
$\C^3$ defined by the equation
$F_d(x,y)=z^m$, where $F_d\in\C[x,y]$ is a homogeneous
form of degree $d$
without multiple factors. If $d\ge 3$ and
$(d,m)\neq (3,2)$
then any closed rational curve $C\subseteq V$
passes through the singular point ${\bar 0}\in V$.
Consequently, the Diophantine equation $F_d(x,y)=z^m$
has a non-constant
polynomial solution in coprime polynomials if and only if
$d\le 2$ or $(d,m)=(3,2)$.\ecor

\proof Indeed, in that case the surface $V$
has an isolated singularity at the origin,
and the defining equation is weighted homogeneous of degree $d_1=md$
with respect to the weights
$w_1=w_2=m,\,w_3=d$,  so that by \ref{cor ci} $\bar
k(V\backslash\{\bar 0\})=-\infty$ if and only if
$$
N_V=d_1-\sum_{j=1}^{3}w_j=(d-2)m-d<0
\quad \text{if and only if}\quad
d\le 2\quad \text{or}\quad (d,m)= (3,2)\,.
$$
Thus the first part follows from Theorem \ref{lines}
above.

If $d=2$ resp., $(d,m)= (3,2)$ then in appropriate coordinates
the equation of the surface $V$ can be written as
$x^2+y^2+z^m=0$ resp., $x^3+y^3+z^2=0$, i.e.,
$V$ is linearly isomorphic to the dihedral surface
$V_{2,2,m}\cong \C^2/\Z_m$ resp.,
to the tetrahedral surface\footnote{The latter is true because
the automorphism group of the projective line $\pP^1$
acts thrice transitively.} $V_{2,3,3}=\C^2/G$ with the
tetrahedral group
$G\subseteq {\rm GL}(2,\C)$ (isomorphic to the alternating group $A_4$).
Thus, there exist polynomial curves
$\C\to V\backslash \{{\bar 0}\}$ (for concrete examples of such
curves see the  Introduction).
\qed

\brem
Corollary \ref{schm} for $(d,m)\neq (3,3),\,(4,3)$ or $(d,2)$
with $d\le 16$
also follows from Theorem 1 in \cite{Sch}
(cf. \cite[Schmidt's Lemma]{KalZa}).
\erem


\bigskip
\begin{thebibliography}{99}

\bibitem[Ab]{Ab} Shreeram S. Abhyankar, {\em Quasirational
singularities}, Amer. J. Math. {\bf 101} (1979), no. 2, 267--300.


\bibitem[Ara]{Ara} D. Arapura,  {\em  A note on Koll\'ar's theorem}, Duke
Math. J.  {\bf 53} (1986), 1125-1130.

\bibitem[Ar]{Ar} M. Artin, {\em  On isolated rational singularities of
surfaces}, Amer. J. Math.  {\bf 88} (1966), 129--136.

\bibitem[BaDw]{BalDw} F. Baldassarri,  B. Dwork,  {\em On
second order linear differential equations with algebraic
solutions}, Amer. J. Math. {\bf 101} (1979), 42--76.

\bibitem[BaKa]{BarKau} G. Barthel, L. Kaup,
{\em Topologie des surfaces complexes compactes singuli\`eres},
in: {\em Sur la topologie des surfaces complexes compactes},
S\'em. Math. Sup. {\bf 80}, Presses Univ. Montr\'eal,
Montreal, Que.,  1982, 61--297.

\bibitem[BPV]{BPVdV} W. Barth, C. Peters, A. van de Ven,
{\em Compact complex surfaces},  Ergebnisse der Mathematik
und ihrer Grenzgebiete (3), Springer-Verlag, Berlin-New
York, 1984.

\bibitem[Be$_1$]{Be1} J. Bertin,
{\em Automorphismes des surfaces non compl\`etes, groupes Fuchsiens et
singularit\'es quasihomog\`enes}, S\'emin. d'alg\`ebre P. Dubreil et M.-P.
Malliavin, 36\`eme Ann\'ee, Proc., Paris 1983/84,  Lect. Notes Math. {\bf
1146} (1985), 106-126.

\bibitem[Be$_2$]{Be2} J. Bertin,  {\em  Pinceaux de droites et
automorphismes des surfaces affines}, J. Reine Angew. Math. {\bf 341}
(1983), 32--53.

\bibitem[Beu]{Beu} F. Beukers, {\em The Diophantine equation
$Ax\sp p+By\sp q=Cz\sp r$}, Duke Math. J. {\bf 91} (1998), no. 1,
61--88.




\bibitem[Bri]{Bri} E. Brieskorn,  {\em Rationale Singularit\"aten
komplexer Fl\"achen}, Invent. Math. {\bf 4} (1968), 336-358.

\bibitem[BrMa]{BrMa} D. Brownawell, D. W. Masser, {\em Vanishing sums
in function fields}, Math. Proc. Cambridge Philos. Soc. {\bf 100}
(1986), no. 3, 427--434.


\bibitem[Bu]{Bur}
D. Burns,  {\em On rational singularities in dimensions $>2$},
Math. Ann. {\bf 211} (1974), 237--244.


\bibitem[ClGr]{CleGri} C. H. Clemens, Ph. A. Griffiths,
{\em The intermediate Jacobian of  the cubic threefold},
Ann. of Math.
{\bf 95} (1972), 281-356.

\bibitem[DaGi]{DanGi} V.I. Danilov, M.H. Gizatullin,  {\em
Automorphisms of affine surfaces. I, II},  Math. USSR Izv. {\bf 9}
(1975), 493-534; ibid. {\bf 11} (1977), 51-98.

\bibitem[DaGr]{DarGr} H. Darmon, A. Granville,   {\em  On the
equations $z\sp m=F(x,y)$ and $Ax\sp p+By\sp q=Cz\sp r$}, Bull. London
Math. Soc. {\bf 27} (1995), 513--543.




\bibitem[Ev]{Ev} A. Evyatar (formerly A. Gutwirth), {\em On polynomial
equations}, Israel J. Math. {\bf 10} (1971), 321--326.

\bibitem[Fl$_1$]{Fl} H.~Flenner, {\em Rationale
quasihomogene Singularitaeten},  Arch. Math. {\bf 36}
(1981), 35-44.

\bibitem[Fl$_2$]{Fl2} H.\ Flenner, {\it
Divisorenklassengruppen quasihomogener Singularit\"aten,}
Journal reine angew.\  Math.\ {\bf 328}  (1981), 128--160.

\bibitem[FlVo]{FlVo} H.Flenner,  W. Vogel, {\em Joins and
intersections}, Springer Verlag, Berlin--Heidelberg--New York, 1999.

\bibitem[FlZa$_1$]{FlZa1} H.\ Flenner,  M.~Zaidenberg,
{\em Log-canonical forms and log-canonical singularities},
e-print math.AG/0012122, 22p.

\bibitem[FlZa$_2$]{FlZa2} H.\ Flenner,  M.~Zaidenberg,
{\em Locally nilpotent derivations on affine surfaces with
a $\C^*$-action} (in preparation).

\bibitem[FA]{FA} {\em Flips and abundance for algebraic threefolds},
Papers from the Second Summer Seminar on Algebraic Geometry
held at the University of Utah, Salt Lake City, Utah, August 1991.
Ast\'erisque {\bf 211} (1992).


\bibitem[GoWa]{Go} S.\ Goto;  Keiichi Watanabe,
{\em On graded rings. I.} J. Math. Soc. Japan {\bf 30},
(1978) 179-213.

\bibitem[GrRi]{GrRi} H.\ Grauert, O.~Riemenschneider,
{\em Verschwindungss\"atze f\"ur analytische Kohomologiegruppen
auf komplexen R\"aumen}, Invent. Math. {\bf 11} (1970), 263-292.


\bibitem[Gur]{Gur} R.V.~Gurjar, {\em
Remarks on the topology of surface singularities and
applications.} With an appendix by K.~Paranjape.  Proceedings of
the Indo-French Conference on Geometry (Bombay, 1989), 95--102,
Hindustan Book Agency, Delhi, 1993.

\bibitem[Ha]{Ha} G.H.\ Halphen, {\em Sur la reduction des
\'equations diff\'erentielles lin\'eaires aux formes
int\'egrables}, M\'emoires pr\'esent\'es par divers savants \`a
l'Academie des sciences de l'Institut National de France, T.
XXVIII, N 1, Paris, F. Krantz, 1883; Oeuvres. Vol. {\bf 3},
Paris 1921, 1-260.


\bibitem[KaZa]{KalZa} S. Kaliman, M. Zaidenberg, {\em
Miyanishi's characterization of the affine 3-space does not
hold in higher dimensions}, Ann. Inst. Fourier (Grenoble) {\bf 50} (2000),
1649--1669.

\bibitem[Kaw]{Kaw} Y. Kawamata, {\em Crepant blowing-up
of $3$-dimensional canonical singularities
and its application to degenerations of surfaces},
Ann. of Math. (2) {\bf 127} (1988), 93-163.

\bibitem[KKMS]{KKMSD}  G. Kempf,  F. Knudsen, D. Mumford, B.
Saint-Donat,   {\em Toroidal embeddings}, Lecture Notes in Mathematics
{\bf 339}, Springer-Verlag, Berlin-Heidelberg-New York, 1973.

\bibitem[Kl]{Kl} F. Klein, {\em Vorlesungen \"uber das Ikosaeder und
die Aufl\"osung der Gleichungen vom f\"unften Grade}, Teubner, Leipzig,
1884.  English transl.: F. Klein, {\em Lectures on the Icosahedron  and
the solution of equations of fifth degree},  Dover, 1956.


\bibitem[Kol]{Kol} J. Koll\'ar,  {\em  Rational curves on algebraic
varieties}, Ergebnisse der Mathematik {\bf 32}, Springer-Verlag,
Berlin, 1995.

\bibitem[Kol1]{Kol1} J. Koll\'ar,  {\em  Higher direct images of dualizing
sheaves}, I. Ann. of Math. (2) {\bf 123} (1986), 11-42; II. {\em ibid.}
{\bf 124} (1986), 171-202.

\bibitem[La$_1$]{La1} H. B. Laufer, {\em Normal two-dimensional
singularities}, Annals of Mathematics Studies.   {\bf 71} (1971),
Princeton, N. J.: Princeton University Press and University of Tokyo
Press. XI, 161 p.

\bibitem[La$_2$]{La2} H. B. Laufer, {\em On rational singularities},
Amer. J. Math. {\bf 94} (1972), 597-608.

\bibitem[ML]{ML} L. Makar-Limanov, {\em  On the hypersurface $x + x^2y
+ z^2 + t^3 = 0$ in ${\C}^{4}$ or a ${\C}^3$-like threefold which is
not ${\C}^3$},  Israel J. Math. {\bf 96} (1996), 419--429.

\bibitem[MaMi]{MaMiy} K. Masuda, M. Miyanishi,
   {\em \'Etale endomorphisms of algebraic surfaces with $G_m$-actions},
  Math. Ann. {\bf 319} (2001), 493--516.

\bibitem[Mi]{Miy} M. Miyanishi,
   {\em Singularities of normal affine surfaces containing 
cylinderlike open sets}
J. Algebra 68 (1981), 268--275.

\bibitem[Mora]{Mora} M. Morales,  {\em Resolution  of
quasi-homogeneous singularities and plurigenera},  Compos.
Math. {\bf 64} (1987), 311-327.

\bibitem[Moi]{Moi} B.G. Moishezon, {\em
On n-dimensional compact complex varieties
with n algebraically independent
meromorphic functions. I-III},
Amer. Math. Soc. Transl. II. Ser. 63 (1967), 51-177.

\bibitem[Moriw]{Moriw} A. Moriwaki, {\em Torsion freeness of
higher direct images of canonical bundles}, Math. Ann. {\bf
276} (1987), 385-398.

\bibitem[Or]{Or} S. Yu. Orevkov, {\em  On singularities that are
quasirational in the sense of Abhyankar}, Uspekhi Mat. Nauk {\bf 50}
(1995), no. 6(306), 201--202.

\bibitem[OW]{OW} P. Orlik, Ph. Wagreich,
  {\em  Equivariant resolution of singularities with $C\sp{*}$ action},
Proceedings of the Second Conference on Compact Transformation Groups
(Univ. Massachusetts, Amherst, Mass., 1971), Part I, pp. 270--290. 
Lecture Notes
in Math.  {\bf 298},  Springer, Berlin, 1972.


\bibitem[Pi$_1$]{Pi1} H. Pinkham,  {\em  Normal surface singularities with
$\C^*$ action}, Math. Ann. {\bf 227} (1977), 183-193.

\bibitem[Pi$_2$]{Pi2} H. Pinkham,  {\em Singularit\'es de Klein---I, II},
S\'eminaire sur les Singularit\'es des Surfaces (Centre de
Math\'ematiques de l'\'Ecole Polytechnique, Palaiseau, 1976--1977), M.
Demazure,  H. Ch. Pinkham and B. Teissier, eds. Lect. Notes Mathem.
{\bf 777}, Springer, Berlin,  1980, 1--20;   {\em Singularites
rationnelles de surfaces. Appendice}, ibid, 147--178.

\bibitem[Ren]{Ren} R. Rentschler, {\em Op\'erations du groupe additif
sur le plane affine}, C.R. Acad. Sci. Paris, {\bf 267} (1968), 384--387.

\bibitem[Sai]{Sai} M. Saito, {\em Modules de Hodge polarisables},
Publ. Res. Inst. Math. Sci.
{\bf 24} (1988), 849-995 (1989).

\bibitem[Sak]{Sak} F. Sakai,  {\em Kodaira dimensions of complements  of
divisors}, Complex analysis and algebraic geometry,  Iwanami Shoten,
Tokyo, 1977,  239--257.

\bibitem[Sch]{Sch} W. M. Schmidt, {\em Polynomial solutions of
$F(x,\,y)=z^{n}$}, Proceedings of the Queen's Number Theory Conference,
1979 (Kingston, Ont., 1979), Queen's Papers in Pure and Appl. Math.
{\bf 54},  Queen's Univ., Kingston, Ont., 1980, 33--65.

\bibitem[Schw]{Schw} H. A. Schwartz, {\em \"Uber diejenigen F\"alle, in
welchen die Gaussische hypergeometrische Reihe eine algebraische
Function ihres vierten Elementes darstellt}, J. f\"ur die reine und
angewandte Mathematik {\bf 75} (1873), 292-335.

\bibitem[Sei]{Sei} A. Seidenberg,  {\em Derivations and integral closure},
Pacific J. Math.  {\bf 16} (1966), 167--173.

\bibitem[Ste]{Ste} S.A. Stepanov,  {\em Diophantine equations over
function fields}, Mat. Zametki  {\bf 32} (1982), 753-764; English transl.
in: Math. Notes {\bf 32} (1982), 861-868.


\bibitem[TSH]{TSH} T. Tomaru, H. Saito, T. Higuchi,   {\em
Pluri-genera $\delta \sb{m}$ of normal surface singularities with
$\C^*$-action}, Sci. Rep. Yokohama Nat. Univ. Sect. I  {\bf 28}
(1981), 35--43.

\bibitem[Vas]{Vas} W. V. Vasconcelos,  {\em
Derivations of commutative noetherian rings},
Math. Z.  {\bf 112}  (1969), 229-233.

\bibitem[Vie]{Vie}
E. Viehweg, {\em Rational singularities of higher dimensional schemes},
Proc. Amer. Math. Soc. {\bf 63}  (1977), 6--8.

\bibitem[Wah]{Wah} J. M. Wahl,   {\em Equations defining rational
singularities}, Ann. Sci. \'Ecole Norm. Sup. (4)  {\bf 10} (1977),
231-263.

\bibitem[KeWa]{KeWa} Keiichi Watanabe,  {\em Certain invariant subrings
are Gorenstein. I, II}, Osaka J. Math. {\bf 11} (1974), 1-8; ibid.
{\bf 11} (1974), 379-388.

\bibitem[KiWa]{KiWa} Kimio Watanabe,  {\em On plurigenera of normal
isolated singularities}, Math. Ann.  {\bf 250} (1980), 65-94.

\bibitem[Za]{Za}  M. Zaidenberg, {\em Isotrivial families of curves on
affine surfaces and characterization of the affine plane}, 
Math. USSR Izvestiya
{\bf 30} (1988), 503-531. Addendum: ibid, {\bf 38} (1992), 435-437.


\end{thebibliography}
\end{document}